\documentclass{amsart}
\input intervals.sty

\usepackage[english]{babel}
\usepackage{graphicx}
\usepackage{amsthm}
\usepackage{amssymb}
\usepackage{amsmath}
\usepackage[usenames,dvipsnames]{color}

\newcommand{ \Proof }{\noindent\textbf{Proof: }}
\newcommand{ \ax }{\ensuremath{\tilde{x}_0} }           
\newcommand{ \PM }{\ensuremath{\mathcal{P}} }         

 \usepackage{epsfig}
 \usepackage{psfrag}
 \newtheorem{teor}{Theorem}
 
 \newtheorem{lema}{Lemma}
 \newtheorem{corol}{Corollary}

 \newtheorem{remark}{Remark}

\def\Zset{\mathbb{Z}}
\def\Nset{\mathbb{N}}
\def\Rset{\mathbb{R}}
\def\Cset{\mathbb{C}}

\def\Pcal{\mathcal{P}}
\def\Ccal{\mathcal{C}}
\def\ii{\mbox{i\,}}
\def\bq{{\bf{q}}}
\def\bp{{\bf{p}}}
\newcommand{\fin}{\hfill $\Box$}

\title{Rigorous KAM results around arbitrary periodic orbits for Hamiltonian
 Systems.}
\author{Tomasz Kapela}

\address{Tomasz Kapela, Jagiellonian University, Institute of Computer Science,
\mbox{Nawojki 11}, 30-072 Krak\'ow, Poland \\
Uppsala University, Department of Mathematics, Box 480, 75106 Uppsala, Sweden }
\email{kapela@ii.uj.edu.pl}

\author{Carles Sim\'{o}}
\address{Carles Sim\'{o}, Dept. de Matem\`{a}tica Aplicada i An\`{a}lisi, 
 Univ. de Barcelona, Gran Via 585, 08007 Barcelona, Spain}
\email{carles@maia.ub.es}

\begin{document}

\begin{abstract}
We set up a methodology for computer assisted proofs of the existence and the
KAM stability of an arbitrary periodic orbit for Hamiltonian systems. We give
two examples of application for systems with 2 and 3 degrees of freedom. The
first example verifies the existence of tiny elliptic islands inside large
chaotic domains for a quartic potential. In the 3-body problem we prove the KAM
stability of the well-known figure eight orbit and two selected orbits of the
so called family of rotating Eights. Some additional theoretical and numerical
information is also given for the dynamics of both examples.
\end{abstract}


\maketitle

\section{Introduction}

KAM Theorem (see \cite{Kol,Arn,Mos} and also \cite{AA}) is a fundamental result
for Hamiltonian systems because it ensures the existence of a set, nowhere dense
but of positive measure, of points of the phase space which behave in a regular,
quasi-periodic way. The main point is that the system should be a perturbation
of an integrable system and that a non-degeneracy condition, asking for the
invertibility of the  actions to frequencies map, has to be satisfied.
The standard notation, being a Hamiltonian system given by
\begin{equation}
  \label{eq:ham}
  \dot{q} =  H_p \qquad  \dot{p} =  - H_q\,,
\end{equation}
where the Hamiltonian $H(q, p):\Omega \to \Rset$ is a smooth function 
defined on an open set $\Omega \subset \Rset^{n+1} \times \Rset^{n+1}$,
will be used.

If we consider the dynamics close to a fixed point the methodology is simple.
Assume that the fixed point is totally elliptic or the problem can be reduced to
the totally elliptic case, for instance by restricting the attention to the
centre manifold. Then one can proceed to compute the normal form up to a
moderate order, say to order 4 in the $(q,p)$ variables. Assuming that
no resonances occur up to this order then one can consider the normal form as
the integrable Hamiltonian and the remainder as the perturbation, and it is
easy to check the non-degeneracy condition. This approach has been used, e.g.,
in the study of the vicinity of the collinear libration points in the general
planar three-body problem, restricted to the centre manifold, see \cite{MaSa}.
A moderate number of arithmetic operations allows to decide on the applicability
of KAM Theorem.

The problem is much more involved when we want to apply KAM Theorem around an
arbitrary totally elliptic periodic orbit which is not known analytically.
Even if some analytic expression of the orbit is available, the study of the
dynamics on the vicinity at the required order can be not feasible analytically.
As it is usual, one can restrict the problem to the study of the vicinity of a
fixed point of a symplectic map on a suitable Poincar\'e section in dimension
$2n$.

The goal of this paper is to set up a methodology for the rigorous check of the
KAM conditions for the symplectic map (see, e.g., \cite{AA}).

We give two examples of application. The first one is a simple classical
Hamiltonian system with two degrees of freedom and depending on a parameter
$c$. The main feature is that the potential consists only on quartic terms.
Changing $c$ the system can be integrable or display large chaotic domains.
In these domains one can guess, by numerical computation of iterates of a
Poincar\'e map, that some tiny islands exist. The problem is to show,
rigorously, that indeed there are elliptic periodic points of the Poincar\'e 
map inside these islands and KAM conditions hold.

The second example concerns the well-known figure eight solution of the general
three-body problem with equal masses. See \cite{CM} for a proof of the existence
of that orbit, found numerically by Moore \cite{Moore}. This is an example of
``choreography'' (see \cite{SimECM}, where the notion of choreographic solution
was introduced, and the references therein), that is, a $T$-periodic solution of
the $N$-body problem ($N=3$ in the present case) such that all the bodies move 
along the same path with time shift $T/N$ between consecutive bodies. This
topic for $N=3$ has been studied by present authors in \cite{KS} where, in
particular, it was proved the totally elliptic character of the figure eight on
fixed energy levels and remaining at the zero level of angular momentum. Using
reductions the problem becomes a Hamiltonian with three degrees of freedom. The
related Poincar\'e map is 4D. In \cite{SimEva} it was claimed, based on a
non-rigorous high order computation of a normal form, that the KAM condition is
satisfied around the figure eight. The computation of the local expansion of the
Poincar\'e map was done by numerical differentiation using multiple precision
and optimal step size for the different orders. In the present paper the
validity of the KAM condition for the figure eight orbit is established
rigorously. 

The paper is organised as follows. In the next two sections the examples with
2 and 3 degrees of freedom are presented and several relevant properties of them
are proved or mentioned. In particular the reduction of the three-body problem
in present case is explicitly carried out, based on \cite{Whi}. Then the
methodology to be applied is explained, introducing the required notation and
emphasizing the rigorous aspects of the CAP (Computer Assisted Proofs). Finally
the results obtained  by applying the methods to both examples are shown.

\section{ A family of quartic potentials} \label{sec:quartic}

As a first example we consider a very simple Hamiltonian
\[ H=\frac{1}{2}(p_x^2+p_y^2+x^4+cx^2y^2+y^4),\]
where $c>-2$ is a real parameter. This is a system widely considered as a
paradigm of chaotic system for $c$ large in the relations between classical and
quantum mechanics, see for instance \cite{BTU,PBB} and references therein. The
energy should be positive and, due to the homogeneity, it can be considered
equal to a fixed value. We shall consider the level $H=\frac{1}{2}$. Note that
for $c\leq -2$ unbounded motion occurs. The system has some obvious symmetries:
It is reversible with respect to time and the changes of sign of $x$ and/or $y$
leave the equations invariant. Furthermore, the symplectic change induced by the
change of variables $(x,y)\to (u,v)=\frac{1}{\sqrt{2}}(x+y,y-x)$ keeps the form
of the Hamiltonian with the parameter $c$ replaced by $\hat{c}=\displaystyle{
\frac{12-2c}{2+c}}$, after a scaling to normalise the coefficients to $x^4$ and
$y^4$ to 1. Obviously the map $c\to\hat{c}$ is an involution having
$c=2$ as fixed point. It can be written also as $\hat{c}+2=\displaystyle{
\frac{16}{c+2}}$. When $c$ ranges in $(-2,\infty)$ increasing its value, the
parameter $\hat{c}$ ranges in the same interval but decreasing.

It is immediate to check that the planes $(x,p_x)$, $(y,p_y)$, $(x=y,p_x=p_y)$
and $(x=-y,$ $p_x=-p_y)$ are invariant and, for the last two, modulo the change
$c\leftrightarrow \hat{c}$, we have the same phase portrait than for the first
two.

The first question to be addressed is the integrability of the Hamiltonian. To
this end we observe that $H=T+V$ is a classical Hamiltonian and $V$ is 
homogeneous of degree $k=4$. We specialize to this degree of homogeneity a
theorem due to Morales-Ramis \cite{MR}

\begin{teor} Assume $H=T+V=\frac{1}{2}(p_1^2+\ldots+p_n^2)+V(x_1,\ldots,x_n)$
where $V$ is homogeneous of degree 4. Let $z=(z_1,\ldots,z_n)$ be a solution
of $z=\nabla V(z)$ and let $\lambda_1,\ldots,\lambda_n$ be the eigenvalues of
${\mbox{Hess }} V(z)$. Then, if $H$ is completely integrable with meromorphic
first integrals, the values of $\lambda_i$ must be equal to numbers of the
form $2m^2-m, 2m^2+\frac{4}{3}m+\frac{7}{72}$ or $2m^2+2m+\frac{3}{8}$ for
$m\in\Zset$.
\end{teor}

\begin{corol} The Hamiltonian $H$ is only integrable for $c=0,2,6$. \end{corol}

\noindent Proof. Beyond the trivial eigenvalue $\lambda_1=3$, due to the
homogeneity, one has $\lambda_2=c/2$, which should be of one of the forms above.
Also $\hat{c}/2$ should be of one of these forms. This reduces the possible
values of $c$ to $0,2,6$. On the other hand the case $c=0$ is obviously 
separable and it is also $\hat{c}=0$ which corresponds to $c=6$. Finally in the
case $c=\hat{c}=2$ the symplectic change $x=r\sin\varphi,\,y=r\cos\varphi,\,p_x=
p_\varphi\cos\varphi/r+p_r\sin\varphi,\,p_y=-p_\varphi\sin\varphi/r+p_r
\cos\varphi$ converts the Hamiltonian to $H=\frac{1}{2} \left(p_\varphi^2/r^2+
p_r^2+r^4\right)$, reducible to 1 degree of freedom.
\fin

\bigskip

\begin{figure}[ht]
\begin{center}
\begin{tabular}{cc}
\hspace{-5mm}\epsfig{file=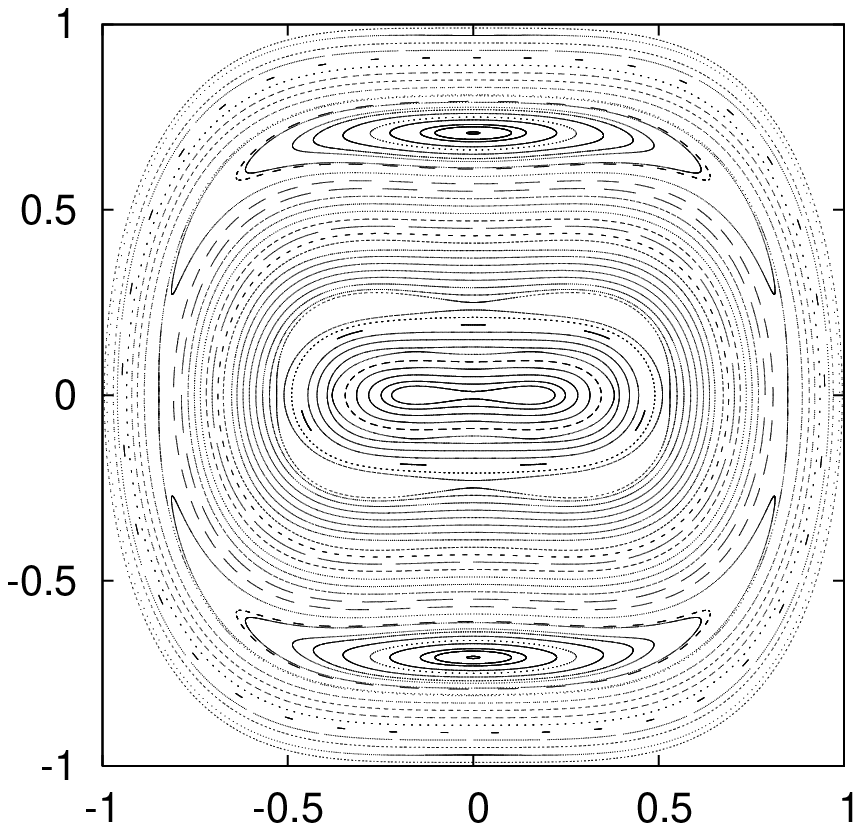,width=7.5cm} & 
\hspace{-26mm}\epsfig{file=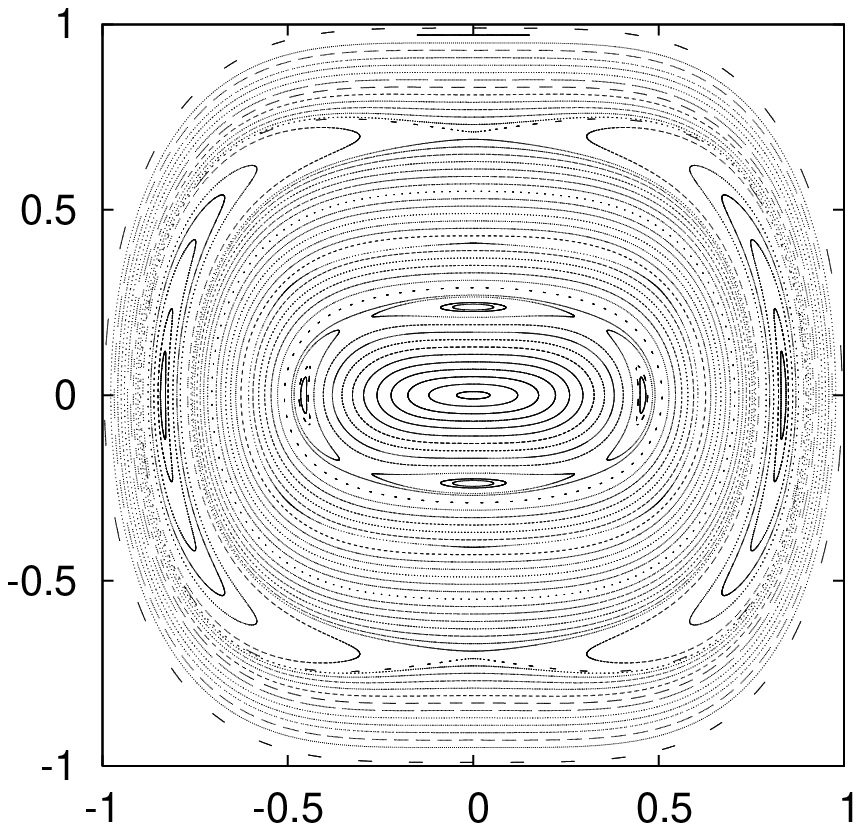,width=7.5cm} \\
\hspace{-5mm}\epsfig{file=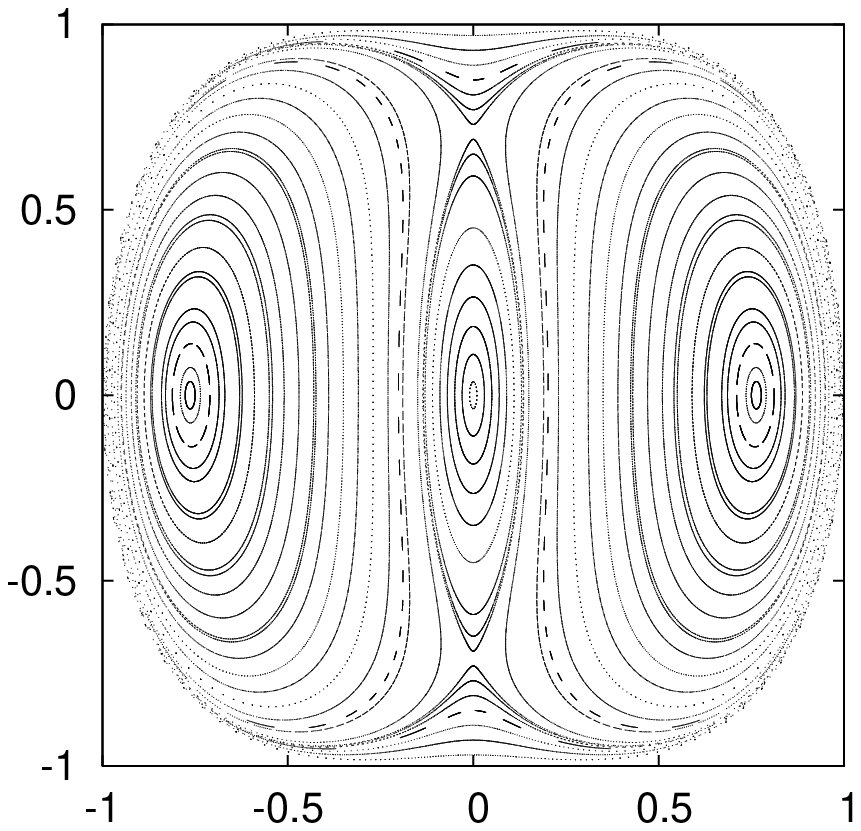,width=7.5cm} &
\hspace{-26mm}\epsfig{file=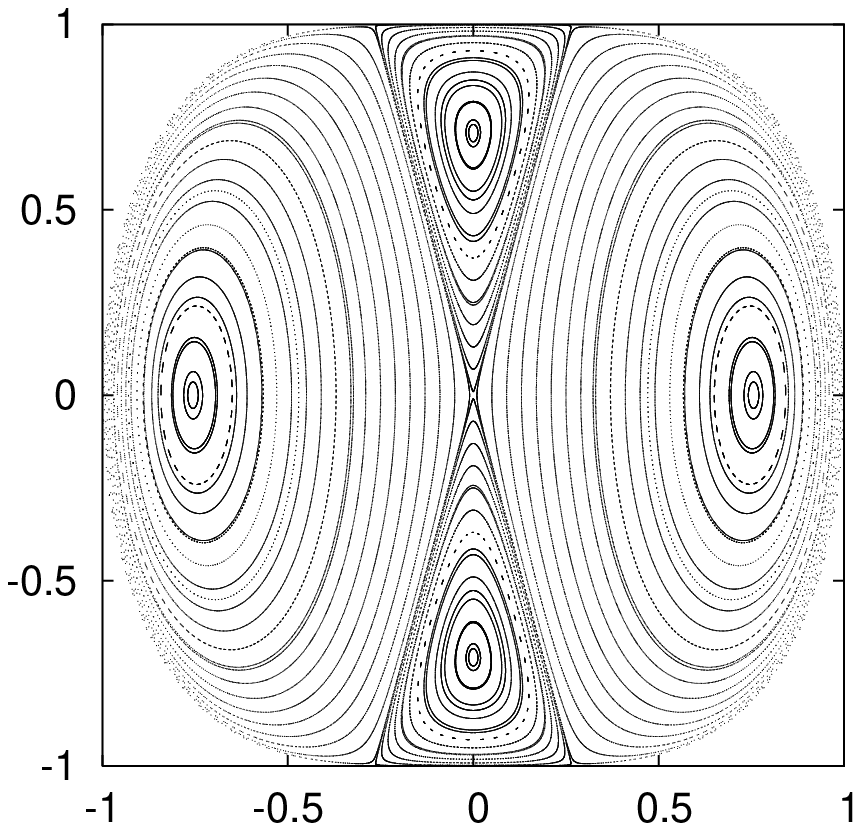,width=7.5cm} \\ 
\hspace{-5mm}\epsfig{file=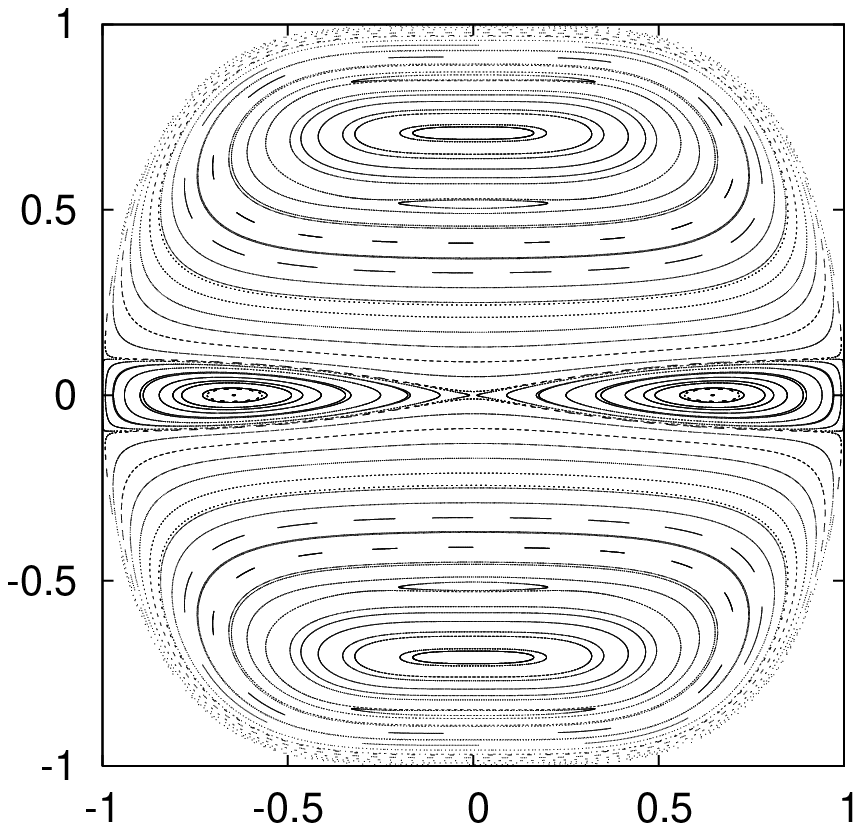,width=7.5cm} &
\hspace{-26mm}\epsfig{file=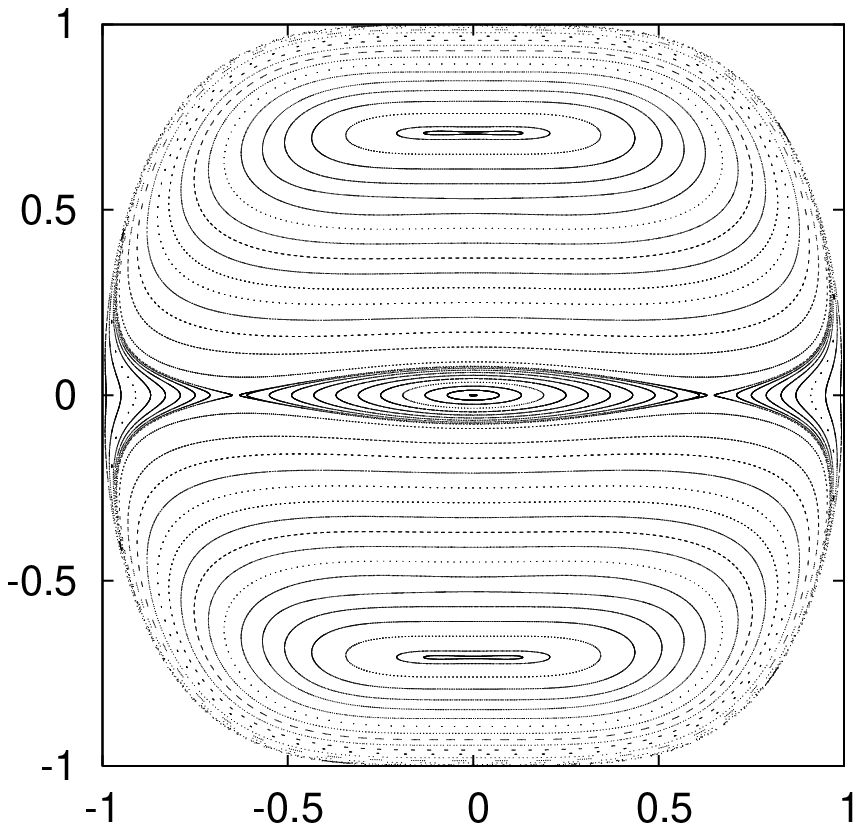,width=7.5cm} 
\end{tabular}
\end{center}
\vspace{-5mm}
\caption{ \label{fig:h4_1} Phase portrait on $\Sigma$, the Poincar\'e section
$y=0$, close to the integrable cases. Left column from top to bottom:
$c=-0.1,1.9,5.9$. Right column $c=0.1,2.1,5.9$. The variables displayed are
$(x,p_x)$.}
\end{figure}

See Figure \ref{fig:h4_1} for an illustration of the phase portrait of $\Pcal$
for values of $c$ close to the integrable cases. In the left (resp. right)
column the value of $c$ is smaller (resp. larger) than the integrable one. The
plots allow to identify easily the bifurcations which occur at the integrable
cases. 

\begin{figure}[htbp]
\begin{center}
\epsfig{file=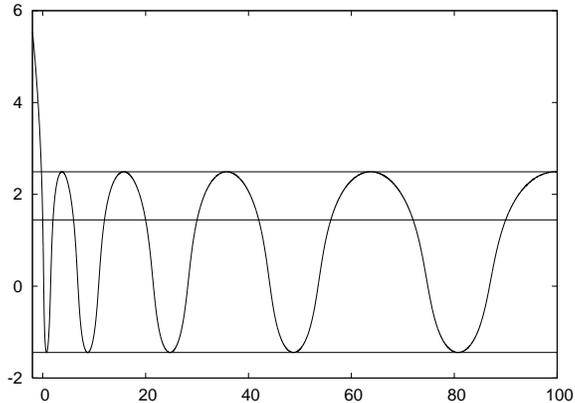,width=8cm}
\end{center}
\vspace{-4mm}
\caption{\label{fig:h4_2}Trace of the differential of the Poincar\'e map $\Pcal$
at the origin as a function of $c$. On the vertical axis, to prevent from the
effect of the large values when $c\to -2^+$, we plot ${\mbox{argsinh}}
({\mbox{Tr}})$ instead of Tr. The horizontal lines correspond to $-2$, $2$ and
$6$.}
\end{figure}

To understand the dynamics of $H$ as a function of $c$ it is useful to consider
the Poincar\'e section $\Sigma$ through $y=0,p_y>0$, defined for $p_x^2+x^4<1$.
The boundary of that domain is a periodic orbit on the invariant plane $y=p_y=
0$. The initial data $x=p_x=0$ give rise to a periodic orbit $y=y(t)$, sitting
on the $(y,p_y)$ plane, which corresponds to a fixed point of the Poincar\'e map
$\Pcal$. The normal variational equations $\xi'=\eta,\eta'=-cy(t)^2\xi$ give the
linear stability of $(0,0)$. It is easy to check that the trace of $D\Pcal$ at
the origin decreases from $+\infty$ to 2 for $c<0$ and then, for increasing $c$,
it oscillates between -2 and 6. In particular it takes the value $+2$ for
$c=c_k=k(k+1), k\in\Nset$. For the values $c=(k-1/2)(k+1/2)$ it takes,
alternatively, the values $+6$ and $-2$ for $k$ even and odd, respectively. See
Figure \ref{fig:h4_2} for an illustration of the behavior of the trace. The set
of values of $c$ for which $|{\mbox{Tr}}|<2$ in the union of the intervals
$(c_0,c_1)\cup(c_2,c_3)\cup \ldots$, where one has linear stability. The first
linear stability intervals are $(0,2),\,(6,12),\,(20,30),\ldots$. It is also
clear that the periodic orbit at the boundary of $\Sigma$ has the same stability
properties as the orbit we have just discussed.

In a similar way one can consider the initial conditions $x=0,\,p_x=\sqrt{1/2}$
which correspond to a periodic orbit in $x=y,p_x=p_y$. In that case one can take
$x+y=0$ as Poincar\'e section. Replacing $c$ by $\hat{c}$ one has a similar
result: For $\hat{c}\to -2^+$, which implies $c\to\infty$, this fixed point is
unstable for $c>6$. Stability intervals for $\hat{c}$ are identical to the ones
given before for $c$ in the case of the fixed point at the origin. They give
rise to intervals which accumulate to $-2$. The first intervals (in the $c$
parameter) are $(2,6),\,(-6/7,0),\,(-3/2,-14/11),\,(-50/29,-18/11),\,
(-42/23),\,(-66/37)$, etc. We shall refer to these two periodic orbits as the
basic ones. Of course, symmetries give rise to similar orbits, like the one 
through $x=0,\,p_x=-\sqrt{1/2}$ 

\begin{figure}[htbp]
\begin{center}
\epsfig{file=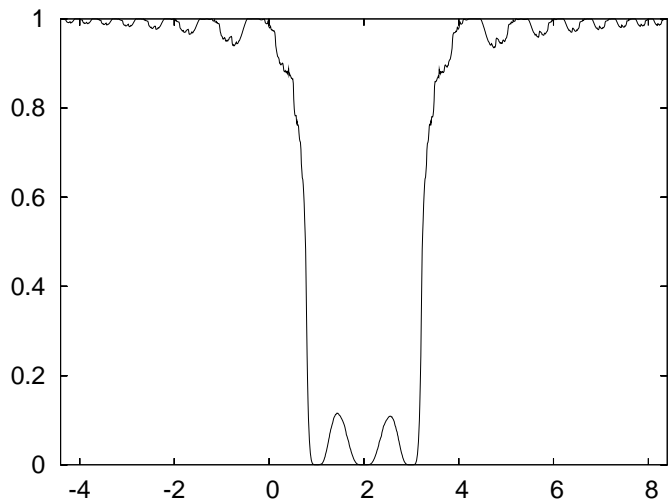,width=10cm}
\end{center}
\vspace{-5mm}
\caption{\label{fig:h4_3}Fraction of chaotic dynamics $\psi(c)$ in the
Poincar\'e section $\Sigma$ as a function of the parameter $c$. To better see
the behavior of the function $\psi$ the variable displayed in the horizontal
axis is $d(c)= \log(c+2)/\log(2)$ instead of $c$. The three integrable cases
$c=0,\,2,\,6$ for which $\psi(c)=0$ are seen on $d=1,\,2$ and $3$,respectively.}
\end{figure}

To see the evolution of the phase space away from the integrable cases we have 
computed an estimate of the ``fraction of chaotic motion'' in $\Sigma$ as a
function of $c$. Due to the symmetries it is enough to do the computations for
$x\ge 0,\,p_x\ge 0$. In the domain bounded by $x=0,\,p_x=0$ and $p_x^2+x^4=1$ we
have selected ``pixels'' with centre of the form $(i/2000,j/2000),i,j\in\Nset$
and for each one we have estimated the maximal Lyapunov exponent $\Lambda$. In
fact we are not interested on the concrete value but rather on whether one can
accept $\Lambda=0$. Symmetry and some other simplifications allow to reduce the
computational task. Then the fraction of chaotic motion $\psi(c)$ is estimated
as the number of pixels for which one has evidence that $\Lambda>0$ divided by
the total number of pixels in the domain. Several checks have been done using
different strategies and maximal number of iterates of the Poincar\'e map to
have reliable information (see, e.g., \cite{PS} for details). The results are
shown in Figure \ref{fig:h4_3}.

The interpretation of the figure is clear: close to the integrable cases most of
the dynamics is regular. In fact, if $c^*$ is one of the values of $c$ for which
one has integrability, the fraction of chaoticity seems to be exponentially
small in $|c-c^*|$ for nearby values. Between these values of $c$ the value of
$\psi(c)$ is below $0.12$. Then, to the left of $c=0$ and to the right of $c=6$
there is a quick increase of $\psi(c)$. But at the ranges in which one of the
basic periodic orbits is linearly stable, we can expect the existence of
islands, which decrease the measure of the chaotic domain. The oscillations of
the decrease are becoming smaller when the limits, either $c\to -2^+$ or $c\to
\infty$ are approached.

\begin{figure}[htbp]
\begin{center}
\begin{tabular}{cc}
\hspace{-10mm}\epsfig{file=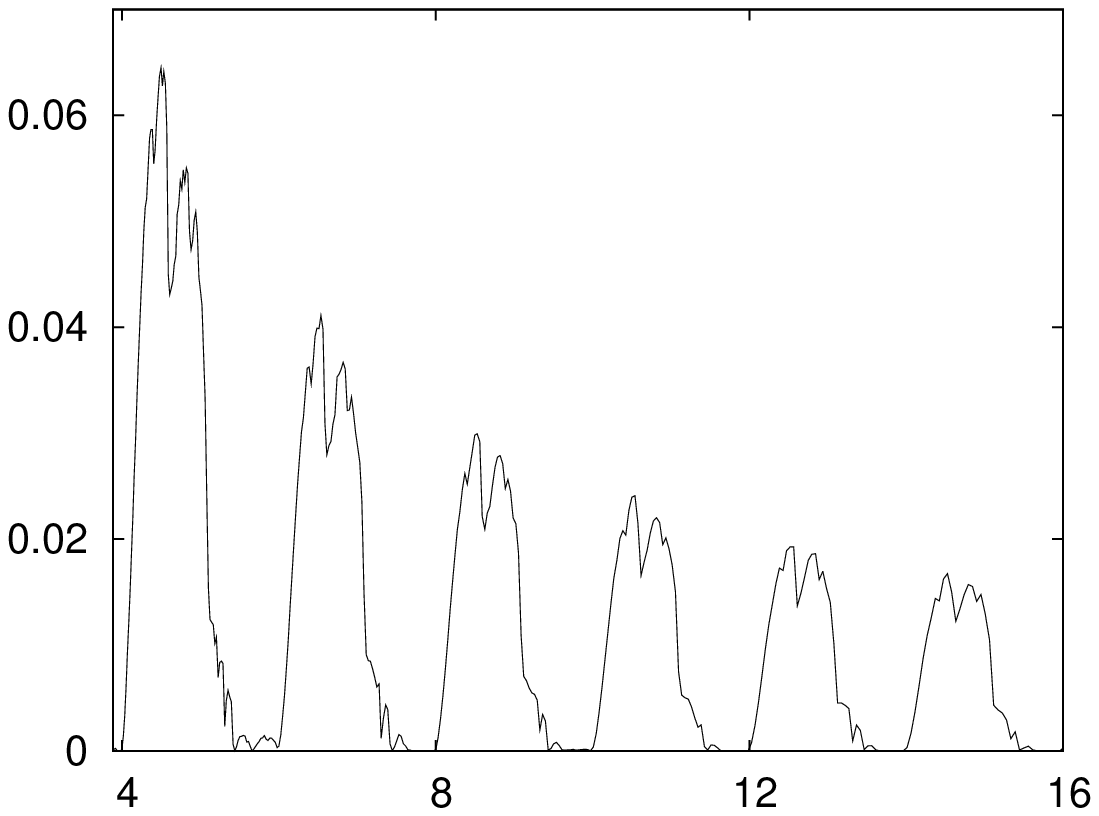,width=8.0cm} & 
\hspace{-10mm}\epsfig{file=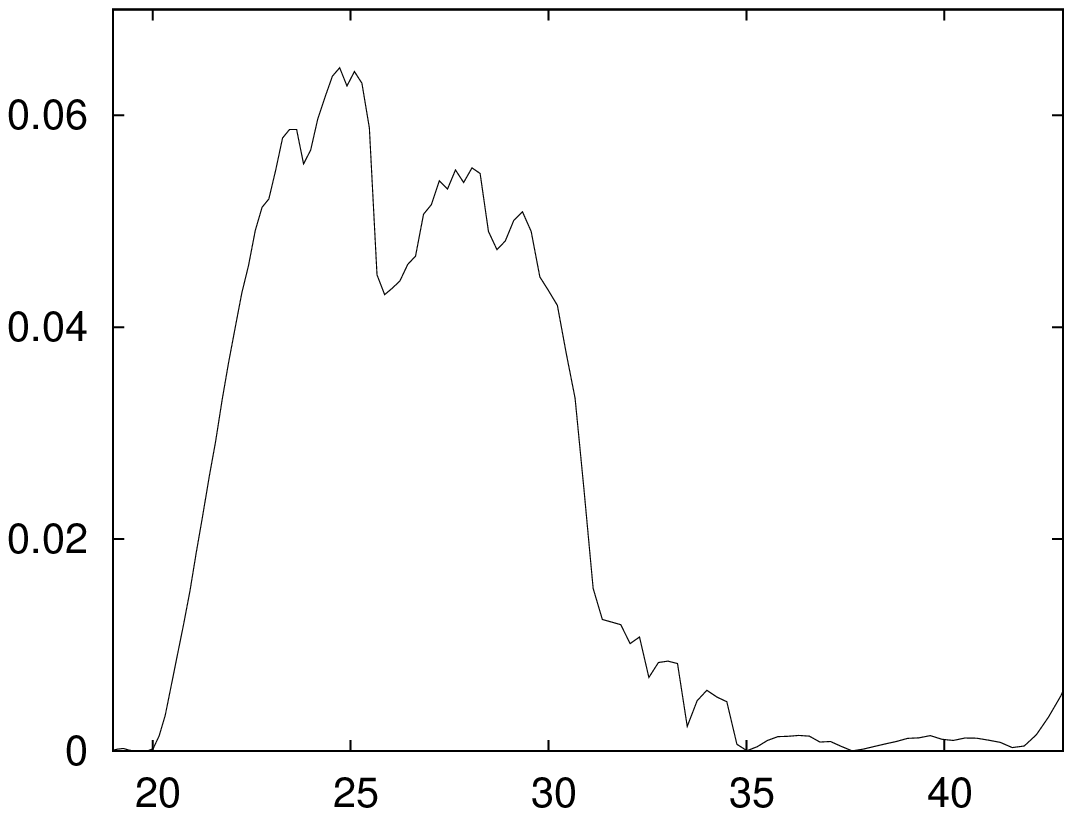,width=8.0cm}  
\end{tabular}
\end{center}
\vspace{-5mm}
\caption{\label{fig:h4_b} Left: The fraction of integrability as a function of
$r$, defined by $r(r+1)=c$. Right: Idem for a range of $c$ containing $[c_4,c_6]
=[20,42]$ using $c$ as horizontal variable.}
\end{figure}

More concretely, the system away from the range of values of $c$ where it is
close to integrable, has a well defined and repetitive structure. To this end
we consider the ``fraction of integrability'' $1-\psi(c)$. Figure \ref{fig:h4_b}
shows it for the values of the parameter $c$ on the right hand side of the
domain which contains the close to integrable dynamics. On the left hand side
of that domain the plot is symmetrical to the one shown in Figure \ref{fig:h4_b}
left. To better see the scaling properties, the horizontal parameter $r$ is an
extension of the index $k$ to the real numbers, defined by $r(r+1)=c$. The
different ``bumps'' are quite similar, and the heights scale as $\approx 1/r$.
The right part of Figure \ref{fig:h4_b} shows the behavior for $c\in[c_4,c_6]$,
which is similar for all the ranges $[c_{2m},c_{2m+2}],m\geq 2.$ At $c=20$ the
central periodic orbit becomes again elliptic (see Figure \ref{fig:h4_2})
(remember that the same thing happens for the periodic orbit at the boundary of
$\Sigma$). A domain of regular motion starts which is reminiscent of the
behavior found in many other problems. See, for instance, the reference
\cite{SV} in the case of the H\'enon map. However there are important
differences due to the many symmetries of present problem. Anyway the mechanisms
to explain the behavior of the plot are the ones explained in \cite{SV}. For
instance, around $c=25.6$ the ``last'' curve surrounding the period 4 islands
breaks down and this increases the size of the main chaotic zone.

\begin{figure}[ht]
\begin{center}
\begin{tabular}{ccc}
\hspace{-5mm}\epsfig{file=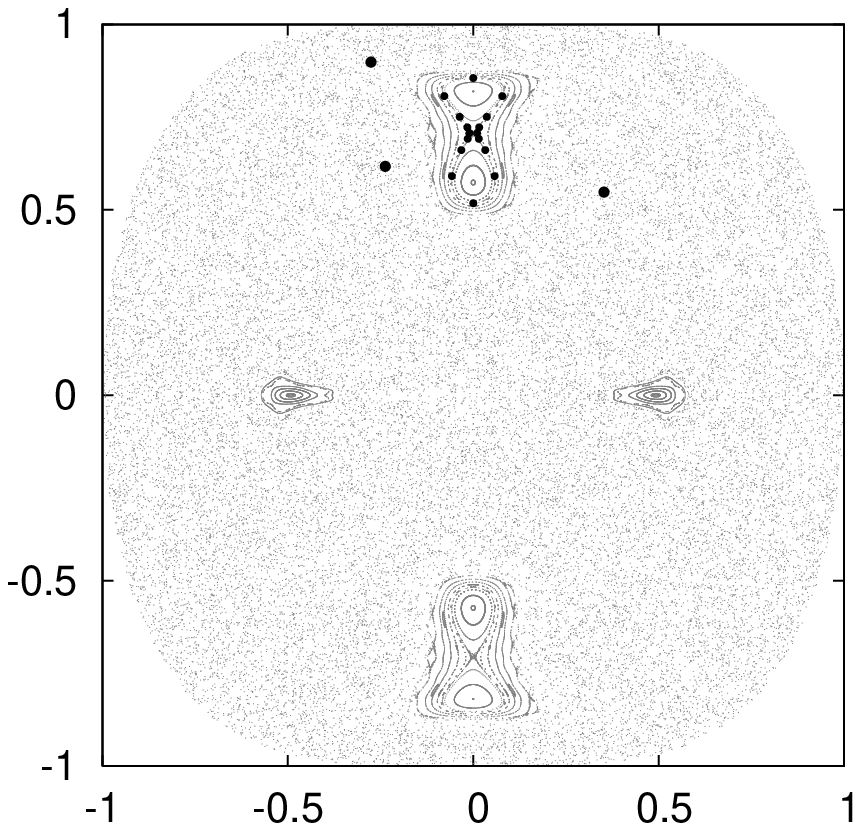,width=7.5cm} & 
\hspace{-26mm}\epsfig{file=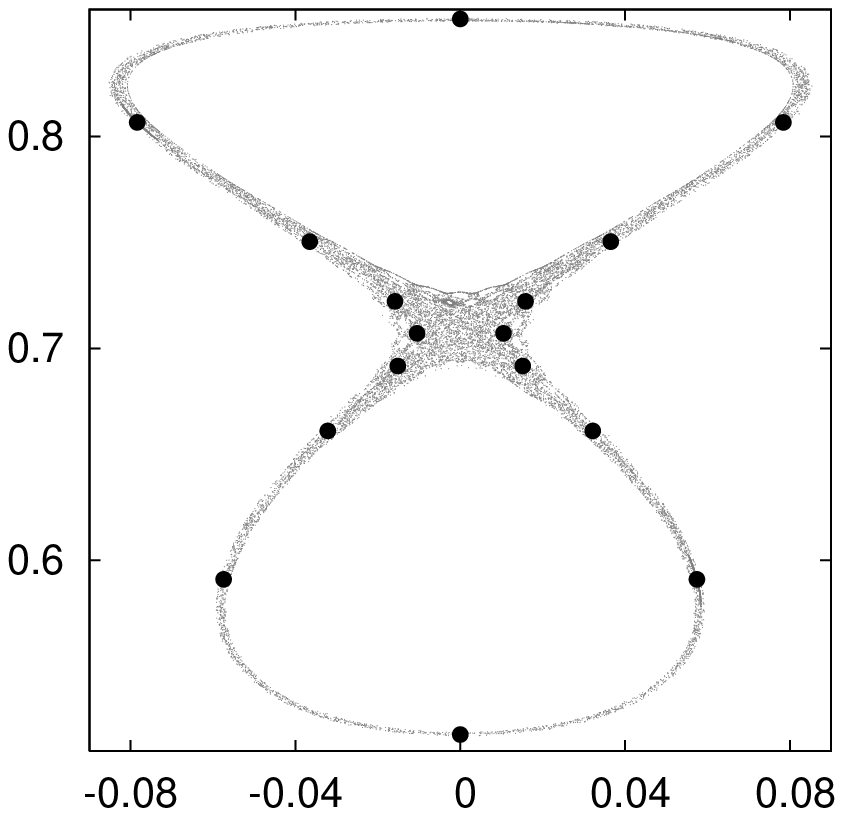,width=7.5cm} & 
\hspace{-26mm}\epsfig{file=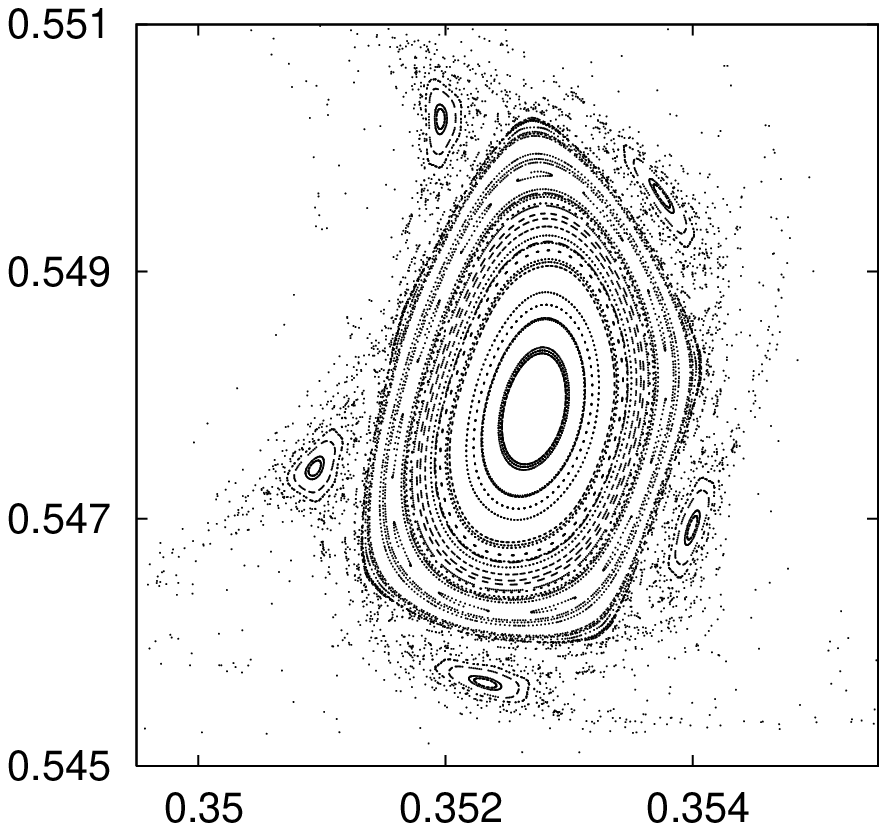,width=7.5cm} 
\end{tabular}
\end{center}
\vspace{-5mm}
\caption{\label{fig:h4_4}Phase portrait on $\Sigma$ for $c=-0.90$. Left: a
global view (in grey) and two periodic orbits (in black). Centre: a
magnification of the chaotic zone around the figure eight separatrix associated
to the unstable periodic orbit on $(0, \sqrt{1/2})$ and a stable periodic orbit
of period 16. Right: The island and some satellite islands near one of the
period 3 points, well inside the large chaotic domain. }
\end{figure}

At $c=30$ the central periodic orbit becomes unstable by a pitchfork bifurcation
and its separatrices display the typical figure eight shape, similar to the
example displayed in the central plot in Figure \ref{fig:h4_4}. Each part of the
figure eight, skipping the invariant curves surrounding the full figure eight
separatrix, is really close to the behavior of the H\'enon map (except by 
nonlinear scalings in variables and parameters). The decrease of the regular
dynamics in Figure \ref{fig:h4_b} right from $c=30$ to $c\approx 31.1$ is due to
the destruction of the invariant curves surrounding the full figure eight 
separatrix. Then, at the approximate $c$ values $32.0,\,32.5,\,33.5,\,35.0$
and $36.8$ the decrease in the fraction is due to the destruction of curves
surrounding period 6, 5, 4 and 3 islands and to the period doubling,
respectively. It is worth to mention that for the case of period 3 almost no
point has regular dynamics. The dynamics in that case, including existence of 
tiny period 3 islands and other properties, is in perfect agreement with the
ones of the H\'enon map as studied in \cite{KSZ}.

Our goal is to detect periodic orbits, well inside a chaotic domain, to which
we can apply the methodology to be explained to test the KAM conditions. To
this end we have selected the value $c=-0.90$ (or, equivalently $\hat{c}=
138/11$), for which the value of $\psi$ is approximately $0.9528$. Iterates of
the Poincar\'e map are easily obtained. Using Taylor method at order 30 and a 
maximal relative truncation error of $10^{-21}$ they are computed at an average
rate of $15,000$ iterates per second for that value of $c$.

The Figure \ref{fig:h4_4} shows a global view of the phase portrait and some
details. As this value of $c$ is in one of the instability domains, the periodic
orbit through $(0,\sqrt{1/2})$ is unstable, but $c$ is not too far from the left
end of one of the stability intervals, at $c=-6/7$. Then the chaotic zone around
the figure eight separatrix is still surrounded by invariant curves. Inside that
chaotic zone there exist several stable periodic orbits. On the central part of
Figure \ref{fig:h4_4} one can see a period 16 orbit. We shall test the KAM
conditions for that orbit. An approximate initial condition on $\Sigma$ is
$x=0,\,p_x=0.51780665799545.$ But also a periodic orbit of period 3, far away
from that separatrix, deeply inside the chaotic domain, has been found. An
approximate initial condition for this orbits is $x=0.352718557335,\,
p_x=0.547882838499$, and the applicability of KAM theorem for this orbit will
also be checked.

We remark that the largest area inside an invariant curve around that period-3
island is of the order of $7\times 10^{-6}$. It is not easy to capture this
periodic orbit, but the previous computation of Lyapunov exponents in a grid
with small stepsize is of great help.

\section{The problem of the KAM stability of the Eight} 
\label{sec:eight}

\begin{figure}[ht]
\begin{center}
\epsfig{file=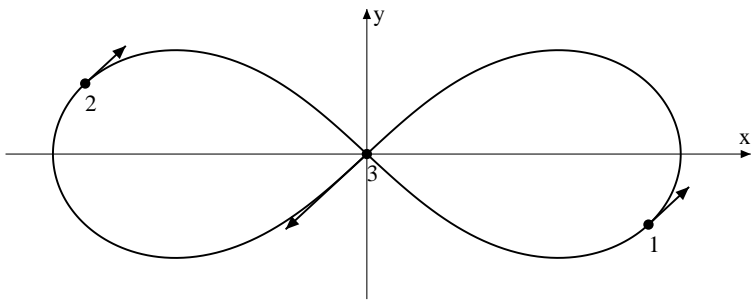,width=10cm}
\end{center}
\vspace{-5mm}
\caption{\label{fig:eight}The figure eight periodic solution for the planar
three-body problem.}
\end{figure}

The figure eight periodic orbit, shown in Figure \ref{fig:eight}, is a
remarkable solution of the planar three-body problem with equal masses
\cite{CM}. The three bodies move on the plane along the same path in solutions
of the form $q(t),q(t+T/3),q(t+2T/3)$ where $T$ is the period. In \cite{SimEva}
a detailed numerical study of this orbit and an extended vicinity of it was
done, looking also at the effect of small changes in the masses and the
bifurcations that they create. See also \cite{CGMS} for choreographies related
to the figure eight, like the satellite and the relative ones. The numerical
evidence given in \cite{SimEva} suggested that non only the orbit was linearly
stable but also KAM theorem applies around it. The rigorous proof of the linear
stability was given in \cite{KS} and the proof of the applicability of KAM is
studied now.

The figure eight solution, three nearby partially hyperbolic periodic orbits,
as well as several 2D and 3D tori around the eight can be seen in \cite{SimEva}.
The orbit has zero total angular momentum and due to the homogeneity of the
potential and Kepler's law, one can fix either the level of (negative) energy
or the period.

The first step to be done is the reduction of the problem to a three degrees of
freedom system. This is a classical result and we follow the exposition that 
can be found in Whittaker's treatise \cite{Whi}. For completeness a short 
account is given below.

\subsection{Reduction of the 3-Body Problem}

We shall assume that the masses are equal to 1. Let $\bq_1,\bq_2,\bq_3$ be the
positions of the three bodies in $\Rset^2$ and $\bp_1,\bp_2,\bp_3$ the
corresponding momenta. The Hamiltonian is 
\[ H=\frac{1}{2}(|\bp_1|^2+|\bp_2|^2+|\bp_3|^2)-\frac{1}{|\bq_2-\bq_3|}-
   \frac{1}{|\bq_3-\bq_1|}-\frac{1}{|\bq_1-\bq_2|}\]
and the angular momentum is ${\bf{M}}=\bq_1\wedge\bp_1+\bq_2\wedge\bp_2+
\bq_3\wedge\bp_3.$

Let $\bq^{\prime}_1,\bq^{\prime}_2,\bq^{\prime}_3,\bp^{\prime}_1,\bp^{\prime}_2
,$ and $\bp^{\prime}_3$ be new variables introduced by means of the generating
function
\[ W=(\bp_1,\bq^{\prime}_1)+(\bp_2,\bq^{\prime}_2)+(\bp_1+\bp_2+\bp_3,
\bq^{\prime}_3),\]
where $(\,\,,\,\,)$ denotes scalar product. We recall that 
\[\bq_j=\frac{\partial W}{\partial \bp_j},\quad \bp^{\prime}_j=\frac{\partial W}
{\partial \bq^{\prime}_j},\; j=1,2,3\] 
and, hence, the change gives
\[ \bq_1=\bq^{\prime}_1+\bq^{\prime}_3,\;\bq_2=\bq^{\prime}_2+\bq^{\prime}_3,\;
\bq_3=\bq^{\prime}_3,\; \bp_1=\bp^{\prime}_1,\; \bp_2=\bp^{\prime}_2,\;
\bp_3=\bp^{\prime}_3-\bp^{\prime}_1-\bp^{\prime}_2.\]
Because of the centre of mass integrals, it is not restrictive to assume $\bq_1+\bq_2+\bq_3=0,\bp_1+\bp_2+\bp_3=0$, which
amounts to $\bq^{\prime}_3=-(\bq^{\prime}_1+\bq^{\prime}_2)/3$ and 
$\bp^{\prime}_3=0$,  Hence, the new
expressions of $H$ and $M$, skipping \mbox{the $^{\prime}$} for simplicity, are
\[ H=|\bp_1|^2+|\bp_2|^2+(\bp_1,\bp_2)-\frac{1}{|\bq_1|}-\frac{1}{|\bq_2|}-
  \frac{1}{|\bq_1-\bq_2|},\quad M= \bq_1\wedge\bp_1+\bq_2\wedge\bp_2,\]
which reduce the system to 4 degrees of freedom.

Let $(q_1,q_2)$ be the components of $\bq_1$, $(q_3,q_4)$ the ones of $\bq_2$
and, in a similar way, we define the components of the $\bp$ variables. We
introduce the generating function
\[ W=
p_1 q^{\prime}_1 \cos q^{\prime}_4+
p_2 q^{\prime}_1 \sin q^{\prime}_4+
p_3 (q^{\prime}_2 \cos q^{\prime}_4-q^{\prime}_3 \sin q^{\prime}_4)+
p_4 (q^{\prime}_2 \sin q^{\prime}_4+q^{\prime}_3 \cos q^{\prime}_4),\]
which gives raise to the transformation
\[ q_1=q^{\prime}_1 \cos q^{\prime}_4,\; 
   q_2=q^{\prime}_1 \sin q^{\prime}_4,\;
   q_3=q^{\prime}_2 \cos q^{\prime}_4-q^{\prime}_3 \sin q^{\prime}_4,\;
   q_4=q^{\prime}_2 \sin q^{\prime}_4+q^{\prime}_3 \cos q^{\prime}_4,\]
\[ p_1=p^{\prime}_1 \cos q^{\prime}_4-((p^{\prime}_4-q^{\prime}_2 p^{\prime}_3+
   q^{\prime}_3 p^{\prime}_2)\sin q^{\prime}_4) /q^{\prime}_1,\]
\[   p_2=p^{\prime}_1 \sin q^{\prime}_4+((p^{\prime}_4-q^{\prime}_2 p^{\prime}_3+
   q^{\prime}_3 p^{\prime}_2)\cos q^{\prime}_4) /q^{\prime}_1,\]
\[ p_3=p^{\prime}_2 \cos q^{\prime}_4-p^{\prime}_3 \sin q^{\prime}_4,\;
   p_4=p^{\prime}_2 \sin q^{\prime}_4+p^{\prime}_3 \cos q^{\prime}_4.\]
Skipping again the $^{\prime}$, one can write the Hamiltonian and angular
momentum in the new variables. It turns out that the new $q_4$ does not appear
in $H$. It is a cyclic variable. Hence, the conjugated variable $p_4$ is
constant. But it is immediate to see that $M=p_4$. We still keep its value in 
the Hamiltonian, despite in the case of the figure eight $M=0$, because it 
plays a role in the ``rotating eights'' solutions. The final reduced Hamiltonian
has the form
\begin{eqnarray} \label{eq:3breduc}
 H&=&p_1^2+p_2^2+p_3^2+p_1p_2-\frac{p_3}{q_1}(q_2p_3-q_3p_2-p_4)+\frac{1}{q_1^2}
   (q_2p_3-q_3p_2-p_4)^2 \\
 &&-\frac{1}{q_1}-\frac{1}{(q_2^2+q_3^2)^{1/2}}-
   \frac{1}{((q_1-q_2)^2+q_3^2)^{1/2}}. \nonumber 
\end{eqnarray}

The new variables have a simple geometrical meaning. Let us denote the positions
of the masses as $m_j,j=1,2,3$ and as $m_im_j$ the vector from $m_i$ to $m_j$.
Then $q_1$ is the norm of $m_3m_1$ and $p_1$ is the component of the linear
momentum of $m_1$ projected along $m_3m_1$; $q_2$ and $q_3$ are the projections
of $m_3m_2$ along $m_3m_1$ and orthogonal to it and, in a similar way, $p_2$ and
$p_3$ are the projections of the linear momentum of $m_2$; finally $q_4$ is the
angle between the $x$-axis and $m_3m_1$. As said, $p_4=M$.

\subsection{Rotating Eights}
\label{sec:roteig}

\begin{figure}[ht]
\begin{center}
\epsfig{file=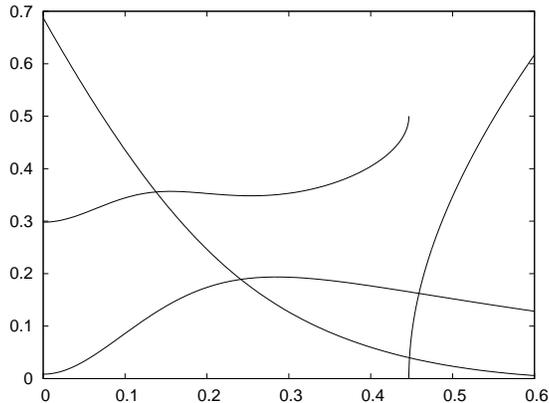,width=8cm}
\end{center}
\vspace{-5mm}
\caption{\label{fig:roteig}Arguments $\alpha_1,\alpha_2$ of the eigenvalues
$\exp(\pm2\pi\ii\alpha_1),\exp(\pm2\pi\ii\alpha_2)$ corresponding to the
rotating Eights solutions for a range of values of angular momentum $M$. When
$\alpha_2$ reaches the value $1/2$ for $M=M_{pd}\approx 0.4467$ a period
doubling is produced and the orbit becomes unstable with a dominant eigenvalue
$\lambda <-1$. From that value on, instead of $\alpha_2$ we plot
$\log(|\lambda|)/2$. The line which decreases monotonically shows the minimal
distance between the bodies along the rotating Eight.}
\end{figure}

If the angular momentum $M$ goes away from zero, the periodic orbit can be 
continued. It becomes quasi-periodic with two basic frequencies, which will
produce a periodic orbit if the ratio of frequencies is rational. But it can be
seen again as a periodic solution, even as a choreography, using a rotating
frame,
a fact noticed by M. H\'enon \cite{Hen} who also found that the continuation
leads to collision orbits. Due to the symmetries of the problem it is enough to
consider $M$ increasing to positive values. For most of these ``rotating
eights'' we can apply the same algorithm as for the Eight. The Figure 
\ref{fig:roteig} shows the values of $\alpha_1<\alpha_2$ such that the
eigenvalues of the Poincar\'e map associated to the rotating eight are
$\exp(\pm2\pi\ii\alpha_1),\;\exp(\pm2\pi\ii\alpha_2),$ as a function of $M$. The
monotonically decreasing line gives the minimal distance between the bodies
along the orbit, an evidence of a nearby collision. The plot changes at $M=
M_{pd}\approx 0.4467$, when $\alpha_2=1/2$, and then the orbit loses stability
with the appearition of a period doubling bifurcation. For $M>M_{pd}$ we plot
the logarithm of the modulus of the dominant eigenvalue (divided by 2 to fit in 
the plot). Low order resonances $k_1\alpha_1+k_2\alpha_2\in\Nset$ with $|k_1|+
|k_2|\leq 4$ are easy to detect for different values of $M$. For them additional
terms will appear in the normal form making our arguments not valid in these
cases. Therefore, as an example, we have selected two values of the angular
momentum, $M_1=0.0048828125$ and $M_2=0.1484375$, that are far enough from
resonances. A non-rigorous exploration of orbits close to the rotating eights
for $|M|<M_{pd}$ gives evidence that for most of them the numerical simulation
suggests the existence of tori, but for some resonances they seem to be
destroyed.

\vspace{-3mm} 

\begin{figure}[ht]
\begin{center}
\begin{tabular}{ccc}
\hspace{-12mm}\epsfig{file=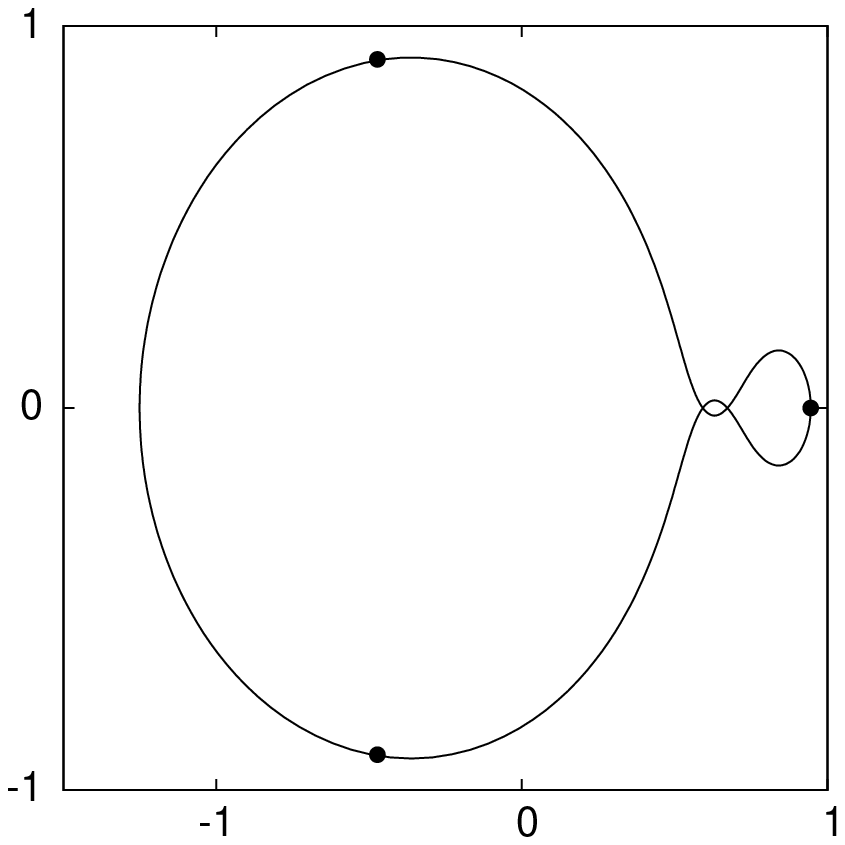,width=6.7cm} & 
\hspace{-20mm}\epsfig{file=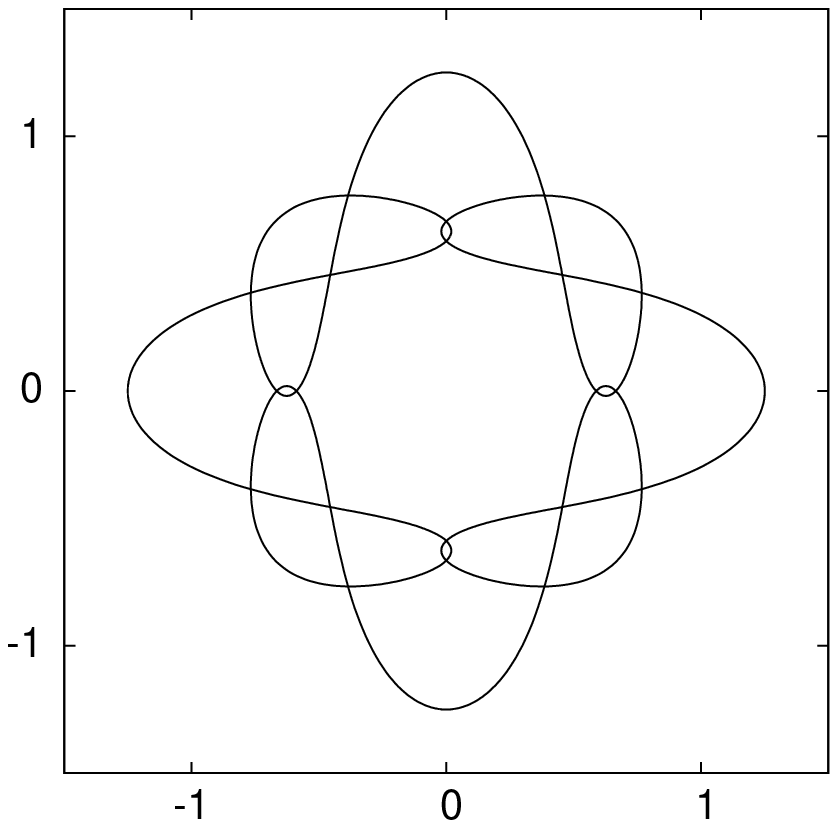,width=6.7cm} &
\hspace{-15mm}\epsfig{file=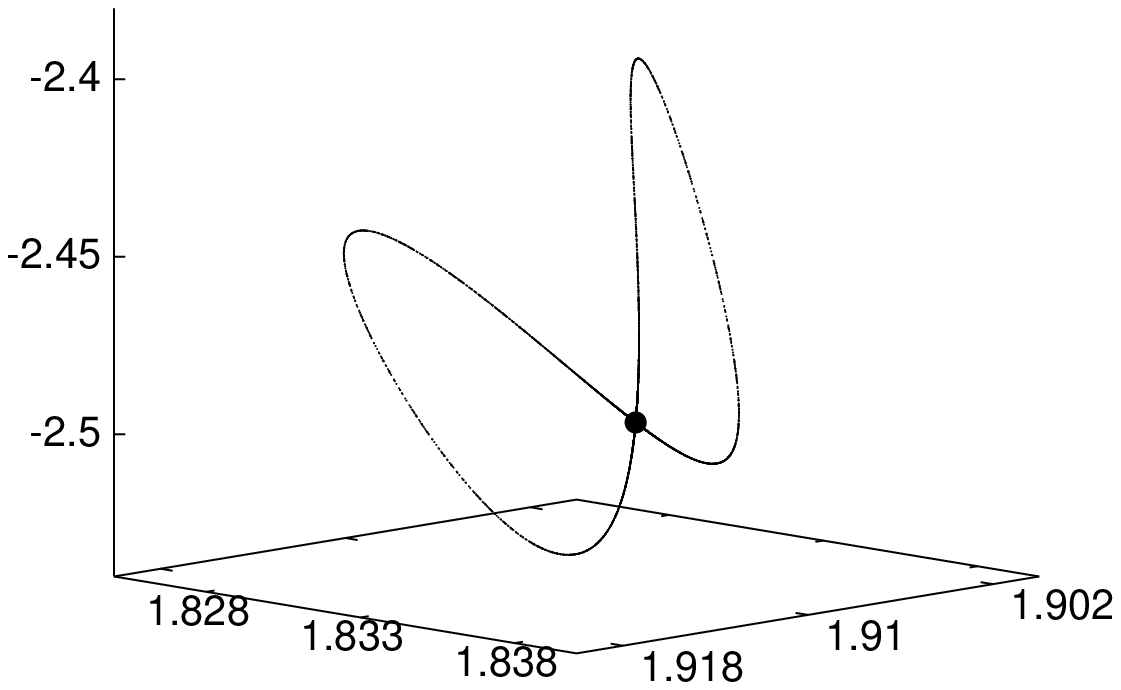,width=6.7cm}
\end{tabular}
\end{center}
\vspace{-5mm}
\caption{\label{fig:orbrota} A rotating eight for $M\approx 0.4493$, already on
the linear instability domain. See the text for a detailed explanation.}
\end{figure}

Initially the rotating eights, in coordinates which rotate with the suitable
angular velocity, look similar to the orbit shown if Figure \ref{fig:eight},
but the left and right hand side lobes are no longer symmetric. Later on the
orbits in the rotating frame can develop extra loops. As an example Figure
\ref{fig:orbrota} shows the solution obtained for $M\approx 0.4493$, shortly
after the orbits become unstable. The value of $M$ has been selected so that in
the fixed frame (middle panel) the orbit is also periodic. On the left panel the
orbit is shown in a rotating frame. At some moment, as displayed, one of the
bodies is on the rightmost point on the path, on the small loop on the right,
and the other two are symmetrical w.r.t. the horizontal axis, on the large loop
on the left. When the bodies move they pass very close to collision on the tiny
central lobe. The right panel shows a 3D projection, on the variables $(q_1,q_2,
p_1),$ of what is seen in the Poincar\'e section. The large dot represents the
periodic orbit. An initial point, taken by adding $10^{-9}$ to $q_1$, is
iterated under the Poincar\'e map. The points are scattered close to the 
``separatrix'' associated to the unstable/stable directions. Of course, these
1D manifolds do not coincide, but the splitting must be expected to be
exponentially small in $M-M_{pd}$.

\section{Preliminaries and the algorithm}
\label{sec:algorithm}

Before going into details we want to emphasize that the algorithm presented in
this paper is \textbf{rigorous}. By rigorous we mean that during all the
computations we take into account and bound all possible errors. In this way we
get not the exact values but verified estimates of the computed quantities.
Therefore the theorems that we apply have assumptions of a special kind (i.e.
inequalities, inclusions etc.) that can be checked using those estimates. As a
result we obtain a computer assisted method to have proofs of the existence and
the KAM stability of periodic orbits for Hamiltonian systems that has full
mathematical rigor. In theory all calculations could be done by paper and
pencil, but in practice the number of operations, even if they are very basic
and trivial, exceeds human resources.

\subsection{Interval arithmetics}

As the precision of the computer is finite we use an interval arithmetic to take
care of the round-off errors. All floating point operations are replaced by the
corresponding operations on closed intervals such that we always obtain some
representable superset of the true result. This is also extended to all
elementary functions.

For the rest of the paper, following \cite{KNN}, we use boldface to denote
intervals and objects with interval coefficients. For those objects we use the
names \emph{interval vector} (or \emph{box}), \emph{interval matrix} etc. to
stress their interval nature. For a set $S$ by $\hull\,S$ we denote the interval
hull of $S$ i.e. the smallest product of intervals containing $S$. For an
interval $\iv{a}=[\lo{a},\hi{a}]$ we define its diameter by $\diam\,\iv{a}:=
\hi{a} - \lo{a}$. We define \emph{the diameter of a box}
$\x = (\x_1, \x_2, \dots, \x_n)$ as a maximal diameter of its components
i.e. $\diam(\x) = \max \{ \diam\, {x}_k | k= 1,\dots,n \}$.

\subsection{The rigorous computation of Taylor expansions of Poincar\'e maps}

To rigorously integrate ODEs and to obtain verified enclosures for the partial
derivatives with respect to the initial conditions  we use $\Ccal^1$-Lohner
\cite{Z_C1} and $\Ccal^r$-Lohner \cite{WZ_Cr} algorithms implemented in the CAPD
library \cite{CAPD}. Those algorithms are based on the Taylor integrator and set
representations proposed by Lohner. For some set of data (say the box $\x$)
and time step $h$  the $\Ccal^r$-Lohner algorithm produces rigorous enclosures
for the solution $\varphi(x,h)$ of the ODE and its derivatives $\frac{\partial^
{|\alpha|} \varphi(x,h)}{\partial x^\alpha}$ for all initial points $x\in\x$ and
all multiindices $|\alpha|\leq r$. The specialized $\Ccal^1$-Lohner algorithm is
able to compute only first order derivatives.

The ``naive'' representation of the derivatives as a table of interval vectors
leads to a huge overestimation due to the ``wrapping effect'' (see
\cite{Mo66,Loh}). Hence, internally during the integration, for each multiindex
$\alpha$ the corresponding vector $v= \frac{\partial^{|\alpha|} \varphi(x,h)}
{\partial x^\alpha}$  is stored in one of the Lohner representations \cite{Loh}.
For that reason even $\Ccal^0$ algorithms need $\Ccal^1$ information to properly
set coordinates to suppress the "wrapping" error. In the current implementation
to store derivatives we use doubletons
      \[v = x_0 + C r_0  + B r,\]
where $x_0$ is a point vector (the centre), $C$ and $B$ are matrices of ``good``
coordinates (usually $C$ is close to the Jacobian matrix and $B$ is its
orthogonalization), $r_0$ represents the initial size of vector $v$ and $r$
stores all computational errors.

To compute all the derivatives up to order $r$ for an $n$ dimensional ODE we
need to solve $n \binom{n+r}{n}$ equations. If they would be solved directly it
would lead to the integration in a high dimensional space and is usually
inefficient (most of the rigorous solvers internally need $C^1$ information that
squares the dimension). The $\Ccal^r$-Lohner algorithm makes use of the special
structure of the variational equations to avoid this and as a result it can
bound derivatives to an arbitrary order $r$ in an efficient way.

On top of those rigorous ODE solvers the CAPD library implements algorithms to
rigorously compute Poincar\'{e} maps and their derivatives with respect to the
initial conditions for an affine sections (for details see \cite{Z_C1, WZ_Cr}).

\subsection{Notation and definitions related to Hamiltonian systems}

We denote by $J$ the Poisson matrix
\[
 J = \left(\begin{array}{cc}
        0 & I_n \\
       -I_n & 0 \\
     \end{array}\right)\,,
\]
where $I_n$ denote the $n$ dimensional identity matrix.
The \emph{Poisson bracket} of functions $f,g:\Rset^{2n} \supset M \mapsto \Rset$
is a new function $\{f, g\} := (\nabla f)^T J \nabla g$.

By $\Nset^+$ we denote set of all positive integers i.e. $\Nset^+=\{1,2,\dots\}$
and we define $\Nset := \{0\} \cup \Nset^+$. An element of $\Nset^n$ will be
called a \emph{multiindex}. For a multiindex $k=(k_1,k_2,\dots,k_n)$ and a
vector $v = (v_1,v_2,\dots,v_n) $ we define
\begin{itemize}
 \item $|k| = \sum_{j=1}^n |k_j|$,
  \item $v^k = v_1^{k_1}v_2^{k_2}\dots v_{n}^{k_{n}}$
\end{itemize}

A vector $(\lambda_1,\lambda_2, \dots, \lambda_n) \in \Cset^n$ satisfies the 
non-resonant condition up to order $r$ if for all multiindices $k=(k_1,k_2,
\dots,k_n)$ such that $1<|k| \leq r$ and all $j\in \{1,2,...,n\}$ we have
\begin{equation}
\label{eq:nonresonant}
 \lambda_j \neq \lambda^k.
\end{equation}

\subsection{The algorithm}

For a given Hamiltonian system \eqref{eq:ham} we  assume that we have an
approximate initial condition $x_0$ of a periodic orbit and a Poincar\'e section
$\Sigma$ that contains $x_0$. For fixed level of energy $h \approx H(x_0)$ the
Poincar\'e map defines a symplectic map
\[
\PM: \Sigma \mapsto \Sigma.
\]
We will use $x=(q,p)\in\Rset^{n}\times \Rset^{n}$ to denote local canonically
conjugated variables.  Let \ax correspond to $x_0$ in these variables.

In this setting to prove existence and KAM stability of a periodic orbit given
by approximate initial conditions $x_0$ it is enough to prove for some $k\in
\Nset^+$ the existence and KAM stability of the unique fixed point of $\PM^k$ in
some small neighbourhood of \ax. The algorithm consists of  the following steps:
\begin{enumerate}
 \item Proof of the existence of the unique fixed point (a periodic orbit).\\
       Rigorous estimates of initial conditions.
 \item Proof of the linear stability. 
 \item Computation of a rigorous Birkhoff normal form.
 \item Checking an appropriate non-degeneracy condition. 
\end{enumerate}
The details for each step will be given in the following subsections. 

\subsection{Proof of the existence}

In the first step of the algorithm we prove the existence of an unique periodic
orbit close to ${x}_0$ and obtain rigorous bounds for its initial conditions.
Therefore the preliminary step is to reduce the Hamiltonian system
\eqref{eq:ham} by suitable symplectic transformations so that for a given energy
level the periodic orbit is isolated. Section \ref{sec:eight} contains a (very
classical) example showing how it was done in the case of the 3-body problem.

For a proof  we take a box $\x \subset \Rset^{2n}$ with centre in $\ax$, and
compute rigorous estimates of the interval Newton operator
       \[ N(\ax, \x, F) = \ax - \hull(DF(\x))^{-1}F(\ax)  \]
for $F(x) = \PM^k(x) - x$. If we succeed to verify that $N(\ax,\x,F)\subset\x$
then the interval Newton theorem \cite{Alefeld94,Neumaier90} ensures that inside
\x there exists a unique $k$-periodic point of \PM\!\!. Instead of Newton method
one can also use the interval Krawczyk method \cite{K,KZ} which do not requires
the whole interval matrix $\hull(DF(\x))$ to be invertible.

\begin{remark}
The problem of proving the existence of zeros of $F$ when, as in present case,
involves the computation of $\PM^k(\x)$ is well suited for the use of the
parallel shooting method. In our implementation we make use of it to improve
precision and to speed up computations.
\end{remark}

\begin{remark}
One can improve rigorous estimates of the initial condition of the periodic
point by further iteration of the interval Newton or Krawczyk operator.
\end{remark}

\subsection{Proof of linear stability}

Current step goal is to prove that all the eigenvalues
$(\lambda_1,\lambda_2,\dots,\lambda_{2n})$ of the differential of the iterated
Poincar\'e map $DP = \frac{\partial\PM^k(\hat{x})}{\partial x}$, where $\hat{x}$
is a $k$ periodic point of $\PM$, lie on the unit circle. We want also to obtain
rigorous estimates $\iv{\lambda_j}$ such that $\lambda_j \in \iv{\lambda_j}$ for
$j=1,\dots,2n$.

The point $\hat{x}$ is not known exactly, but from the previous step we have
rigorous estimates \iv{x} of it. From estimates $\iv{DP}=
\frac{\partial\PM^k(\x)}{\partial x}$ using e.g. verified root finding methods
applied to the characteristic polynomial one can obtain estimates of the
eigenvalues. Because those estimates are given by some boxes, part of them
are out of the unit circle. But still the proof of linear stability is possible
due to the fact that our system is Hamiltonian.
\begin{lema}
Let $A\in \Rset^{2n \times 2n}$ be a symplectic matrix with eigenvalues
$(\lambda_1, \dots, \lambda_{2n})$, and let $\iv{\lambda_j}$ be boxes such that $\lambda_j \in \iv{\lambda_j}$ for $j=1,\dots, 2n$. If the following holds
\begin{itemize}
 \item[(A1)] $0 \notin \im(\iv{\lambda_j}) $ for $j=1,\dots, 2n$
 \item[(A2)] $(\iv{\lambda_j}^{-1} \cap \iv{\bar\lambda_k} \neq \emptyset)
 \implies j = k$,
\end{itemize}
then all eigenvalues of $A$ are distinct and lie on the unit circle.
\end{lema}
\Proof
The matrix $A$ is symplectic, hence  if $\lambda$ is an eigenvalue of $A$ then
also $\lambda_j = \lambda^{-1}$ and $\lambda_k = \bar{\lambda}$ are eigenvalues.
But then assumptions (A1) and (A2) ensure that $\lambda^{-1}=\bar{\lambda}$ and
hence $\|\lambda\|=1$. From (A2) we have also that all eigenvalues are distinct.
\qed

This general method requires sharp bounds for the eigenvalues.
This can be a not so easy task in general. Another possibility is
to translate constraints to the characteristic polynomial. In this case one
proves first that the eigenvalues are on the unit circle and then rigorous
enclosures $\iv{\lambda_i}$ are obtained using this fact (see for example
\cite{KS}).

\subsection{Computation of Birkhoff normal form}

The literature on how to compute Birkhoff normal form is very rich. There are
also general software packages that can do it in a non rigorous way. Here we
want to explain how to make this process rigorous. All the computations are done
using interval arithmetics in the way that at the end the normal form will have
interval coefficients that contain the exact values.

Let $T(x) = (\sum_{k\in \Nset^{2n}} c_{1,k} x^k, \dots, \sum_{k\in \Nset^{2n}}
c_{2n,k} x^k)$ be the Taylor series of an analytic symplectic map around a
totally elliptic fixed point and let $(\lambda_1, \dots, \lambda_n,
\bar{\lambda}_1, \dots, \bar{\lambda}_n)$ be the eigenvalues of the linear part
of $T$. Then for $j\in \{1,2,...,2n\}$ and a multiindex $k$ such that $k_j=k_{j
\pm n}+1$ and $k_m = k_{m+n}$ for $m\neq j$ the condition \eqref{eq:nonresonant}
is not satisfied and a resonance occurs. We call it an \emph{unavoidable
resonance} and we say that the term $c_{j,k} x^k$ corresponds to that resonance.
A Taylor series $T$ is said to be \emph{non-resonant} if only unavoidable
resonances are present.

The goal of this section is to make a symplectic change of variables such that
in the new variables the Taylor series $T$ up to a given order $r$ reads
\[
  z \mapsto \Lambda z +  T_r(z),
\]
where $\Lambda$ is a diagonal matrix, $\Lambda=diag(\lambda_1, \dots, \lambda_n,
\bar{\lambda}_1, \dots, \bar{\lambda}_n)$, and $T_r(z)$ contains only terms
corresponding to unavoidable resonances. This is the so called
\emph{non-resonant Birkhoff normal form}. In what follows we present an
algorithm that computes the non-resonant Birkhoff normal form up to order 3,
but it can be easily extended to any given order. 

First, using $C^r$-Lohner algorithm \cite{WZ_Cr} we compute rigorous enclosures
of the coefficients of the Taylor expansion of $\PM^k$ up to order 3 for all
points in \x (an estimate of the fixed point). As a result we obtain a
symplectic map (up to order 3)
\[T:(q_1,\dots,q_n,p_1,\dots,p_n)\rightarrow (\hat{q_1},\dots,\hat{q_n},
\hat{p_1},\dots,\hat{p_n}).\]
We recall that from the previous step we know that the eigenvalues $\lambda =(
\lambda_1,\dots,\lambda_{2n})$ of the linear part of $T$ lie on the unit circle.
  
As a second step we pass the linear part to a diagonal form
\[\Lambda=\text{diag}(\lambda _1,\lambda_2,\dots,\lambda_{2n}).\]
To this end we use the linear change, $B$, given by a matrix formed by
eigenvectors corresponding to the above eigenvalues. As we do not know the exact
eigenvalues, the estimate  of an eigenvector corresponding to $\lambda_i$, has
to be valid for all $\lambda_i \in \iv{\lambda_i}$. To ensure that $B$ is
symplectic we use the following lemma and replace the previous eigenvectors by
suitable multiples of them.

\begin{lema}
Let $A\in \Rset^{2n \times 2n}$ be a symplectic matrix with eigenvalues 
$(\lambda_1,..., \lambda_n, \lambda_{n+1}=\bar{\lambda}_1, \dots, \lambda_{2n} =
\bar{\lambda}_n)$ such that $|\lambda_i| = 1$ and $\lambda_i \neq \lambda_j$ for
$i,j = 1, \dots, 2n$ and $i \neq j$. Let $(e_1,\dots, e_{2n})$ be corresponding
eigenvectors. If $e_i^T J e_{i+n} = 1$ for $i=1, \dots, n$ then the eigenbasis
$B=(e_1,e_2, \dots, e_{2n})$ is symplectic.
\end{lema}
\Proof
To be symplectic matrix $B$ needs to satisfy $B^T J B = \{e_i^T J e_j \} = J$. 
Due to antisymmetry we can assume that $i\leq j$. For $j = i+n$ we already have
that $e_i^T J e_{i+n} =1$. For $j \neq i+n$ we have $\lambda_i \lambda_j \neq 1$
and therefore $e_i^T J e_j=e_i^T A^T J A e_j=\lambda_i \lambda_j e_i^T J e_j=0.$
\qed

In our implementation initially $e_i$ and $e_{i+n}$ are complex conjugate 
vectors, therefore $e_i^T Je_{i+n}=\ii c$ for some $c\in\Rset$. We want to scale
those vectors to get $e_i^T Je_{i+n}=1$ and additionally $e_i=\ii\bar{e}_{i+n}$.
This is possible only if $c< 0$. Therefore if $c>0$ we simply interchange the
indices of corresponding eigenvalues and eigenvectors. Finally we set
\[
 e_{i}\leftarrow\tfrac{e_{i}}{\sqrt{-c}},\qquad e_{i+n}\leftarrow\tfrac{e_{i+n}}
{-\ii \sqrt{-c}}.
\]

Let $u=(u_1,\dots,u_{2n})$ be the new coordinates in this basis. The above
scaling implies that $u_i = \ii \bar{u}_{i+n}$ for $i=1,\dots,n$. The final form
of the symplectic transformation up to order 3, before starting the normal form
computation, is 
\[S=B^{-1}\circ T\circ B : u \rightarrow \Lambda u + S_2 + S_3,
\]
where $S_2$ and $S_3$ denote quadratic and cubic terms of $S(u)$ respectively. 

Let us denote as $U=(U_1,\dots,U_{2n})$ the coordinates of the normal form. We
achieve a normal form in two steps, by cancelling first all the terms of degree
two in $S$ and then the terms of degree three, except the unavoidable
resonances. For the first step we should select, in principle, a transformation
of the form 
\[u=N_1(U)=U+Q(U),\]
where $Q(U)=(Q_1,\dots,Q_{2n})^T$ are quadratic terms. To cancel the terms of
degree two in $S$ we require that
\begin{equation}
\label{eq:SN1}
S\circ N_1=N_1\circ\Lambda
\end{equation}
holds up to degree two. Let us express this in coordinates. Assume that the
quadratic terms of $S$ and $N_1$ are written, respectively, as
\[
(S_2)_j = \sum_{|k|=2}c_{j,k} u^{k}, \qquad Q_j = \sum_{|k|=2}d_{j,k} U^{k},
\qquad \mbox{for } j=1,...,2n,
\]
where $k$ is a multiindex. The condition (\ref{eq:SN1}) allows to obtain
\[
d_{i,k}=\frac{c_{i,k}}{\lambda_i-\lambda^{k}}
\]
for all the required indices $i, k$. If there are no resonances of order 2, 
then all denominators are different from zero.
We ensure this using rigorous estimates $\iv{\lambda}_i$.

However, the map $N_1$, as it was introduced, is not a symplectic map. Its
differential satisfies the symplecticity conditions only to order 1 in $U$ and
we need additional cubic terms to satisfy them to order 2. This suggests to
define the symplectic transformation as the time-1 map of some Hamiltonian $W$
\[ u=N(U)=\varphi_{t=1}^W(U).\]
To determine $W$ we require
\[ \frac{\partial W}{\partial U_{j+n}}=Q_j,
\quad -\frac{\partial W}{\partial U_j}=Q_{j+n},\quad \mbox{ for } j = 1,\dots,n
\]
As the difference between the map and the identity starts with quadratic terms,
the Hamiltonian starts with cubic terms. Then the time-1 map adds terms of
degree 3 to the initial ones. Finally we obtain the components of $N(U)$ as
\[ (N(U))|_j = U_j + Q_j + \frac{1}{2} \{Q_j, W \},\]
where $\{\;,\;\}$ denotes the Poisson bracket. This produces the normal form
to order two as
\[S:   U\mapsto \Lambda U + S_3(U),\]
where $S_3(U)$ are the cubic terms.

The last step is to remove all cubic terms except those corresponding to
unavoidable resonances. We know that there exists a symplectic, near the
identity transformation that will cancel the non-resonant cubic terms leaving
the resonant terms unchanged. Hence we can simply set to zero all terms in
$S_3(U)$ for which we are able to verify non-resonant condition.

Finally we obtain the normal form to order three (we use $z$ as a new
variables)
\[S:   z\mapsto \Lambda z + T_3(z),\]
where $T_3(z)$ are the cubic terms corresponding to unavoidable resonances.

It remains to check the non-degeneracy condition which allows to apply KAM
theorem. This is now easy and will be done directly on the examples.

\section{Results}

In this section we present the results of the application of the algorithm to
the examples introduced in sections \ref{sec:quartic} and \ref{sec:eight}.
The computed values are often very thin intervals with lower and upper bound
having many common leading digits. To increase readability when printing
intervals we put first those common digits and then the remaining digits of
lower and upper bound as subscript and superscript correspondingly, e.g.
\begin{eqnarray*}
 123.456789_{1234}^{5678} & = & [123.4567891234, 123.4567895678], \\
 -123.456789_{1234}^{5678} & = & [-123.4567895678, -123.4567891234], \\
  {}_{-0.0000000001}^{\;\;\;0.0000000001} & = & [-0.0000000001, 0.0000000001].
\end{eqnarray*}
To be rigorous the presented intervals are rounded outwards, e.g the interval
$[10^{-21},10^{-18}]$ when rounded to 3 decimal places reads $0.00_0^1$. That
sometimes can suggests that consecutive interval operations are not correct but
in fact they are performed with much higher precision than the one displayed,
e.g. the result of $1000 \cdot 0.00_0^1$ can still be equal to $0.00_0^1$.

\subsection{The quartic potential example}

For an approximated initial condition as given  in section \ref{sec:quartic} we
carry out the algorithm described in section \ref{sec:algorithm}. The computed
values are displayed in the Table \ref{tab:h4}.
From the first step of the algorithm we obtain a rigorous estimate of the
initial condition: \x. We use them to get enclosures for the eigenvalues:
$\iv{\lambda_1}$ and $\iv{\lambda_2}$. The Birkhoff normal form of the third and
the sixteenth iterate of the Poincar\'e map, respectively, computed at the fixed
point $x \in \x$ in both cases is proved to be
\begin{equation}
 \left(\begin{array}{l}
z_1\\z_2
\end{array}
\right) \mapsto \left(
\begin{array}{l}
 \lambda_1 z_1 + c_1 z_1^2 z_2  \\
 \lambda_2 z_2 + c_2 z_1 z_2^2 \\
\end{array}
\right)
\end{equation}
for some $\lambda_j \in \iv{\lambda_j}$ and $c_j \in \iv{c_j}$ for $j=1,2$.
It is enough to work only with the first equation, because $z_1=\ii\bar{z}_2$,.
Hence the map (up to order 3) reads
\begin{equation}
\label{ex:h4_twist}
z_1  \mapsto  z_1\exp(2\pi\ii(\alpha+d|z_1|^2))
\end{equation}
where $\lambda_1= \exp(2\pi\mbox{i}\alpha)$ and $d = \frac{c_1}{2\pi\lambda_1}$.
Finally, for the twist map \eqref{ex:h4_twist} if $d$ (the torsion, or twist
coefficient in that case) is different from zero then the fixed point is KAM
stable (see \cite{AA, Mos}). The computed rigorous bounds $\iv{d}$, shown in
Table \ref{tab:h4}, verify that $d$ is not zero for both orbits. Therefore they
are KAM stable.

The proof of the existence of the period 16 orbit using double precision did not
succeed, and we were forced to use multiprecision interval arithmetic and Taylor
method of higher order. This increases in a significant way the computational
time.

\begin{table}
\begin{center}
\begin{tabular}{|c|c|c|}\hline
&& \\
period & 3 & 16 \\
&& \\
precision & double (52 bits) & multiprecision (100 bits)  \\
&& \\
$\x$ &
$\left(\begin{array}{c}
0.35271855733_{27137}^{67754}\\
0.547882838_{4949869}^{5020712}
\end{array}\right)$
&
$\left(\begin{array}{c}
{}_{-0.0000000000000001}^{+0.0000000000000001}\\
0.517806657995457_{3}^{4}
\end{array}\right)$ \\
&& \\
$\diam(\x) $                     & $7 \cdot 10^{-12}$                                         & $10^{-23} $ \\
&& \\
$\iv{\lambda_1},\iv{\lambda_2}$  & $0.539_{4}^{6} + \ii 0.84198_{4}^{7}$                      &   $-0.93239050_4^5 \pm \ii 0.36145255_2^3$\\
&& \\
$\iv{c_1}$                       & $89_{149.071}^{315.580} - \ii 139_{180.278}^{345.393}$    &  $137393_{0.985}^{1.069} - \ii 532621._{066}^{177}$ \\
&& \\
$\iv{c_2}$                       & $-89_{156.837}^{307.807} - \ii 139_{187.959}^{337.719}$   &  $-137393_{0.988}^{1.065} - \ii 532621._{070}^{173}$ \\
&& \\
$\iv{d} $                        & $263_{05.594}^{42.269} + \ii {}_{-18.281}^{+18.283}$  &  $-234523.9_{51}^{71} + \ii {}_{-0.011}^{+0.011}$\\
&& \\
\hline
\end{tabular}
\vspace{3mm}
\caption{ Rigorous estimates of the initial condition  {\boldmath $x$,
eigenvalues $\lambda_1,\lambda_2$, coefficients of normal form $c_1, c_2$ and
torsion $d$}
for the proof of the KAM stability of orbits of period 3 and 16 for the quartic potential.  \label{tab:h4}}
\end{center}
\end{table}

\subsection{The figure eight orbit}

We start with the coordinate system introduced in section \ref{sec:eight} in
which the Hamiltonian of the planar 3-body problem has the form
\eqref{eq:3breduc}. The variables are $(q_1,q_2,q_3,p_1,p_2,p_3)$. For the
figure eight orbit and nearby orbits we select as Poincar\'e section the passage
through $q_3=0$, that is, when the three bodies are aligned. We recall that this
collinear passage happens 6 times in a full revolution. Then, given $q_1,q_2,
p_1$ and $p_2$ and the value of the energy $H=h$ we can determine $p_3>0$.
The Poincar\'e map defines a symplectic 4D map on the fixed level of energy and
$M=0$. Local variables which are canonically conjugate are $q_1,q_2,p_1,p_2$.

Rough initial conditions in these variables if we start when the bodies are
aligned, on the $x$ axis, with $m_1$ to the right and $m_2$ at the center, are
\[ q_1=1.9909837697297968,\, q_2= 0.9954918848648984,\, q_3=0,\,\] 
\[p_1=-0.34790196497952825,\,p_2=0.69580392995905650,\,p_3=1.0678596267584018.\]
Of course, $p_3$ can be recovered from the level of energy, which was fixed to
\[
h=-1.2929708570.
\]
The corresponding period is close to $2\pi$. Of course, by the homogeneity of
the potential, any period or any negative value of $h$ are equivalent. 


Now everything is prepared to carry out the algorithm described in section
\ref{sec:algorithm}. The computations are done using multiprecision intervals
with 200 bits of mantissa. For ODE integration we used Taylor method of order
50. As a result we obtained very sharp rigorous enclosures for the initial
conditions of the Eight
\[ \x = (1.990983769917805_{2}^{3},\;0.995491884958902_{6}^{7},
\;-0.34790196496310_{19}^{20},\;0.69580392992620_{39}^{40})\]
Next we computed the Taylor expansion of the Poincar\'e map at the fixed point
$\hat{x}\in\x$ and we proved that all eigenvalues of the linear part lie on the
unit circle. Their rigorous estimates are:
\[
\iv{\lambda}_1,\iv{\lambda}_3 = 
0.998599982092038_{3}^{4} \mp \ii 0.052896840792060_{6}^{7}
\]
\[
\iv{\lambda}_2,\iv{\lambda}_4 =
-0.29759666751871_{29}^{30} \pm \ii 0.954691690275848_{4}^{5}
\]

The final normal form is

\begin{equation}
 \left(\begin{array}{l}
z_1\\z_2\\z_3\\z_4
\end{array}
\right) \mapsto \left(
\begin{array}{l}
 {\lambda_1} z_1 + c_{11} z_1^2 z_3 +  c_{12} z_1 z_2 z_4 \\
 {\lambda_2} z_2  +  c_{21} z_1 z_2 z_3 + c_{22} z_2^2 z_4\\
 {\lambda_3} z_4 + c_{31} z_1 z_3^2 +  c_{32} z_2 z_3 z_4 \\
 {\lambda_4} z_4  + c_{42} z_1 z_3 z_4 + c_{41} z_2 z_4^2\\
\end{array}
\right)
\end{equation}
where
   \[c_{11} \in -235.97_{4}^{5} + \ii 12._{499}^{500}, \qquad
     c_{12} \in 22.02_{3}^{4} - \ii 1.16_{6}^{7}, \]
   \[c_{21} \in -6.56_{3}^{4} + \ii 21.05_{5}^{6}, \qquad
     c_{21} \in -0.15_{1}^{2} + \ii 0.48_{6}^{7}, \]
   \[c_{31} \in 235.97_{4}^{5} + \ii 12._{499}^{500}, \qquad
     c_{32} \in -22.02_{3}^{4} - \ii 1.16_{6}^{7}, \]
   \[c_{42} \in 6.56_{3}^{4} + \ii 21.05_{5}^{6}, \qquad
     c_{41} \in 0.15_{1}^{2} + \ii 0.48_{6}^{7}\,. \]

Due to the form of the eigenvalues and to the fact that no resonances appear to
order 3, except the unavoidable ones, and using the fact that the second and
fourth variables are related to the complex conjugates of the first and third
ones, we only need to work with the two complex variables $z_1, z_2$.

At this point the map reads as
\begin{equation}
\label{ex:eightNF}
 \left(\begin{array}{l}z_1\\z_2\end{array}\right) \mapsto
\left(\begin{array}{l}\lambda_1 z_1+c_{11}|z_1|^2z_1+c_{12}|z_2|^2z_1\\
         \lambda_2 z_2+c_{21}|z_1|^2z_2+c_{22}|z_2|^2z_2 \end{array}\right).
\end{equation}

Let us write $\lambda_j=\exp(i2\pi\alpha_j)$ for $j=1,2$. Then the map
can be written in the form
\begin{equation}
\label{ex:eight_twist}
\left(\begin{array}{l}z_1\\z_2\end{array}\right)  \mapsto
\left(\begin{array}{l}
  z_1\exp(i2\pi(\alpha_1+d_{11}|z_1|^2+d_{12}|z_2|^2))\\
  z_2 \exp(i2\pi(\alpha_2+d_{21}|z_1|^2+d_{22}|z_2|^2))
\end{array}\right), 
\end{equation}
where $d_{jk}=c_{jk} / \lambda_{j}$ for $j=1,2,\;k=1,2$. In the above expression
we let aside, as in all the computations, terms of order 4 or higher. To obtain
\eqref{ex:eight_twist} we take first $\lambda_1 z_1$ and $\lambda_2 z_2$ as
factors in \eqref{ex:eightNF} and compute logarithms. Note that the values of
the coefficients $d_{jk},\,j=1,2,\,k=1,2 $ must be real, although from rigorous
computations we get complex intervals around some real point. It should also
hold that the matrix formed by the $d_{jk}$ coefficients is symmetric, a fact
which is compatible with the results obtained in the computations.

For the figure eight orbit we have
\[
  d_{11} \in -37.6_{09}^{10} + \ii {}_{-0.001}^{+0.001}, \qquad
  d_{12} \in 3.51_{0}^{1} + \ii {}_{-0.001}^{+0.001},
\]
\[
  d_{21} \in 3.51_{0}^{1} + \ii {}_{-0.001}^{+0.001}, \qquad
  d_{22} \in 0.08_{1}^{2} + \ii {}_{-0.001}^{+0.001}\,.
\]
For the map \eqref{ex:eight_twist} the non-degeneracy KAM condition is simply
that the determinant of the torsion $d = d_{11}d_{22}-d_{12}d_{21}$ must be
different from zero. Because for the Eight we have 
\[
d \in \iv{d} = -15.37_{2}^{3} + \ii {}_{-0.001}^{+0.001}
\]
this finishes the proof that it is KAM stable on the level of fixed energy and
angular momentum $M=0$.

\subsection{Stability of rotating Eights}

Exactly the same algorithm as for the figure eight orbit can be used for the
rotating Eights introduced in section \ref{sec:roteig}. The obtained Birkhoff
normal form and torsion condition are the same. Therefore, keeping the same
notation as in the previous section, we only list in Table \ref{tab:roteig}
the rigorous estimates that verify KAM stability of rotating Eights for the two
selected values of angular momentum $M_1=0.0048828125$ and $M_2=0.1484375$.

\begin{table}[ht]
\begin{center}
\begin{tabular}{|c|c|c|}\hline
&& \\
 M  & 0.0048828125 & 0.1484375000 \\
&& \\
 \x &
$\left(\begin{array}{c}
1.990931659347265_{5}^{6}\\
1.008912014249292_{0}^{1}\\
-0.357203564790858_{3}^{4}\\
0.695944468385828_{3}^{4}
\end{array}\right)$ &
$\left(\begin{array}{c}
1.951511291544930_{4}^{5}\\
1.375174867045771_{8}^{9}\\
-0.706576995371201_{4}^{5}\\
0.835536267754800_{4}^{5}
\end{array}\right)$ \\
 && \\
$\diam(\x)$ &   $10^{-52}$  & $2 \cdot 10^{-30}$\\
&& \\
$\iv{\lambda_1},\iv{\lambda_3}$ &   $ 0.99851219_6^7 \mp \ii 0.05452883_4^5$ & $ 0.65401844_6^7 \mp \ii 0.756478_{599}^{600}$\\
$\iv{\lambda_2},\iv{\lambda_4}$ &   $-0.29874539_1^2 \pm \ii 0.95433285_1^2$ & $-0.62092254_4^5 \pm \ii 0.78387192_4^5$\\
&& \\
$\iv{D}=(\iv{d_{jk}})$ &
$\left(\begin{array}{cc}
-36.39_{4}^{5} + \ii {}_{-0.001}^{+0.001} & -3.50_{5}^{6} + \ii {}_{-0.001}^{+0.001} \\
-3.50_{5}^{6} + \ii {}_{-0.001}^{+0.001} & 0.08_{8}^{9} + \ii {}_{-0.001}^{+0.001} \\
\end{array}\right)$ &
$\left(\begin{array}{cc}
-0.5_{49}^{50} + \ii {}_{-0.001}^{+0.001} &
3.28_{0}^{1} + \ii {}_{-0.001}^{+0.001}   \\
3.28_{0}^{1} + \ii {}_{-0.001}^{+0.001} &
1.55_{7}^{8} + \ii {}_{-0.001}^{+0.001} \\
\end{array}\right)$ \\
 && \\
 $\iv{d}=det(\iv{D})$ & $ -15.49_{6}^{7} + \ii {}_{-0.001}^{+0.001}$ & $-11.61_{7}^{8} + \ii {}_{-0.001}^{+0.001}$\\
&& \\\hline
\end{tabular}
\vspace{3mm}
\caption{ \label{tab:roteig} Rigorous estimates of the initial condition 
{\boldmath $x$, eigenvalues $\lambda_i$ and torsion $d$} for the proof of KAM
stability of rotating Eights with angular momentum $M_1=0.0048828125$ and
$M_2=0.1484375$.}
\end{center}
\end{table}

\section{Acknowledgements}

The research of T.K. was supported in part by Polish Ministry of Science
and Higher Education through grant N201 024 31/2163. The work of C.S. was
supported by grants MTM2006-05849/Consolider (Spain) and 2009 SGR 67
(Catalonia).

\bibliography{thebib}{}

\begin{thebibliography}{10}

\bibitem{Alefeld94}
G.~Alefeld.
\newblock {Inclusion methods for systems of nonlinear equations--the interval
  Newton method and modifications}.
\newblock {\em Topics in Validated Computations}, pages 7--26, 1994.

\bibitem{Arn}
V.~I. Arnold.
\newblock Proof of {A}. {N}. {K}olmogorov's theorem on the preservation of
  quasi-periodic motions under small perturbations of the {H}amiltonian.
\newblock {\em Usp. Mat. Nauk SSSR}, 18:13--40, 1963.

\bibitem{AA}
V.~I. Arnold and A.~Avez.
\newblock {\em Probl\`emes ergodiques de la m\'ecanique classique}.
\newblock Gauthier-Villars, Paris, 1967.

\bibitem{BTU}
O.~Bohigas, S.~Tomsovic, and D.~Ullmo.
\newblock {Manifestations of classical phase space structures in quantum
  mechanics}.
\newblock {\em Physics Reports}, 223:43--133, 1993.

\bibitem{CAPD}
{CAPD}.
\newblock {{CAPD} - {C}omputer {A}ssisted {P}roofs in {D}ynamics, a package for
  rigorous numerics}.

\bibitem{CGMS}
A.~Chenciner, J.~Gerver, R.~Montgomery, and C.~Sim\'o.
\newblock {Simple Choreographic Motions of $N$ Bodies: A Preliminary Study}.
\newblock {\em Geometry, Mechanics and Dynamics, Springer-Verlag}, pages
  287--308, 2002.

\bibitem{CM}
A.~Chenciner and R.~Montgomery.
\newblock A remarkable periodic solution of the three body problem in the case
  of equal masses.
\newblock {\em Annals of Mathematics}, 152:881--901, 2000.

\bibitem{PBB}
G.~G. de~Polavieja, F.~Borondo, and R.~M. Benito.
\newblock {Scars in Groups of Eigenstates in a Classically Chaotic Sistem}.
\newblock {\em Physical Review Letters}, 73:1613--1616, 1994.

\bibitem{Hen}
M.~H\'enon.
\newblock Private communication.
\newblock 2000.

\bibitem{KS}
T.~Kapela and C.~Sim{\'o}.
\newblock {Computer assisted proofs for nonsymmetric planar choreographies and
  for stability of the {E}ight}.
\newblock {\em Nonlinearity}, 20(5):1241--1255, 2007.

\bibitem{KSZ}
T.~Kapela, C.~Sim\'o, and P.~Zgliczy{\'n}ski.
\newblock {Some properties of the H\'enon map in the 3:1 resonance}.
\newblock {\em In preparation}, 2011.

\bibitem{KZ}
T.~Kapela and P.~Zgliczy{\'n}ski.
\newblock {The existence of simple choreographies for the {$N$}-body
  problem---a computer-assisted proof}.
\newblock {\em Nonlinearity}, 16(6):1899--1918, 2003.

\bibitem{KNN}
R.B. Kearfott, M.T. Nakao, A.~Neumaier, S.P. Shary, and P.~van Hentenryck.
\newblock {Standardized notation in interval analysis}.
\newblock {\em Proc. XIII Baikal International School-seminar "Optimization
  methods and their applications", Irkutsk, Baikal, July 2-8, 2005}, 4, 2005.

\bibitem{Kol}
A.~N. Kolmogorov.
\newblock On the conservation of conditionally periodic motions under small
  perturbations of the {H}amiltonian.
\newblock {\em Dokl. Akad. Nauk SSSR}, 98:527--530, 1954.

\bibitem{K}
R.~Krawczyk.
\newblock {Newton-{A}lgorithmen zur {B}estimmung von {N}ullstellen mit
  {F}ehlerschranken}.
\newblock {\em Computing (Arch. Elektron. Rechnen)}, 4:187--201, 1969.

\bibitem{Loh}
R.~J. Lohner.
\newblock Einschliessung der lösung gewonhnlicher anfangs- and
  randwertaufgaben und anwendungen.
\newblock {\em Universität Karlruhe (TH) these}, 1988.

\bibitem{MaSa}
R.~Mart\'{\i}nez and A.~Sam\`a.
\newblock On the centre manifold of collinear points in the planar three-body
  problem.
\newblock {\em Celestial Mechanics \& Dynamical Astronomy}, 85:311--340, 2003.

\bibitem{Moore}
C.~Moore.
\newblock {Braids in Classical Gravity}.
\newblock {\em Physical Review Letters}, 70:3675--3679, 1993.

\bibitem{Mo66}
R.~E. Moore.
\newblock {\em Interval analysis}.
\newblock Prentice-Hall Inc., Englewood Cliffs, N.J., 1966.

\bibitem{MR}
J.J. Morales and J.P. Ramis.
\newblock {A note on the non-integrability of some Hamiltonian systems with a
  homogeneous potential}.
\newblock {\em Methods and Applications of Analysis}, 8:113--120, 2001.

\bibitem{Mos}
J.~K. Moser.
\newblock On invariant curves of area-preserving mappings of an annulus.
\newblock {\em Nachr. Akad. Wiss. G\"ottingen, math.-phys. Kl.}, 11a:1--20,
  1962.

\bibitem{Neumaier90}
A.~Neumaier.
\newblock {\em {Interval methods for systems of equations}}.
\newblock Cambridge Univ. Press, 1990.

\bibitem{PS}
J.~Puig and C.~Sim\'o.
\newblock {Resonance tongues in the Quasi-Periodic Hill-Schr\"odinger Equation
  with three frequencies}.
\newblock {\em Regular and Chaotic Dynamics}, 16:62--79, 2011.

\bibitem{SimEva}
C.~Sim\'o.
\newblock Dynamical properties of the figure eight solution of the three-body
  problem.
\newblock {\em Contemporary Mathematics AMS}, 292:209--228, 2000.

\bibitem{SimECM}
C.~Sim\'o.
\newblock {New families of Solutions in $N$--Body Problems}.
\newblock {\em Proc. 3rd ECM, 2000, Progress in Mathematics series,
  Birk\"auser}, 210:101--115, 2001.

\bibitem{SV}
C.~Sim\'o and A.~Vieiro.
\newblock {Resonant zones, inner and outer splittings in generic and low order
  resonances of Area Preserving Maps}.
\newblock {\em Nonlinearity}, 22:1191--1245, 2009.

\bibitem{Whi}
E.~T. Whittaker.
\newblock {\em A treatise on the Analytical Dynamics of Particles and Rigid
  Bodies}.
\newblock Cambridge Univ. Press, 1970, fourth edition, reprinted.

\bibitem{WZ_Cr}
D.~Wilczak and P.~Zgliczy{\'n}ski.
\newblock {$C^r$-Lohner algorithm}.
\newblock {\em Schedae Informaticae}, 20:9--46, 2011.

\bibitem{Z_C1}
P.~Zgliczynski.
\newblock {$C^1$-Lohner algorithm}.
\newblock {\em Foundations of Computational Mathematics}, 2:429--465, 2008.

\end{thebibliography}
\bibliographystyle{plain}

\end{document}